\documentclass[11pt]{article}%
\usepackage{amssymb}
\usepackage{amsmath}
\usepackage{amsfonts}
\usepackage{graphicx}%
\providecommand{\U}[1]{\protect\rule{.1in}{.1in}}
\begin{document}

\title{\textsf{Axiomatic Justification in Constructive Morse Set Theory}}
\author{\textsf{Douglas S. Bridges}}
\maketitle

\begin{abstract}%
%TCIMACRO{\TeXButton{noindent}{\noindent}}%
%BeginExpansion
\noindent
%EndExpansion
\textsf{Working within Constructive Morse Set Theory (CMST)
\cite{BAlps,dsb2022}, we introduce axioms for a new notion, }$\mathsf{jst}%
\ Pp$\textsf{, intended to capture what it means for }$P$\textsf{ to prove, or
justify, }$p$\textsf{ under the BHK interpretation of intuitionistic logic.
Since it makes no distinction between terms and formulae---every term is also
a formula, and vice versa---CMST is well suited to our axiomatic development
of justification theory within set theory itself. After stating our axioms for
}$\mathsf{jst}\ Pp$\textsf{, we derive many consequences thereof. In
particular, we show that (with certain restrictions) our axioms for
}$\mathsf{jst}\ Pp$\textsf{ align with the intended BHK interpretations of the
axioms of intuitionistic logic.}

\end{abstract}

%

%TCIMACRO{\TeXButton{sf}{\normalfont\sf}}%
%BeginExpansion
\normalfont\sf
%EndExpansion%
%TCIMACRO{\TeXButton{depth 1}{\setcounter{secnumdepth}{0}}}%
%BeginExpansion
\setcounter{secnumdepth}{0}%
%EndExpansion

\subsection{Introduction}%

%TCIMACRO{\TeXButton{noindent}{\noindent}}%
%BeginExpansion
\noindent
%EndExpansion
The context of this paper\footnote{%
%TCIMACRO{\TeXButton{sf}{\normalfont\sf}}%
%BeginExpansion
\normalfont\sf
%EndExpansion
\textbf{Keywords: \ }Morse set theory, constructive, justification
\par
\ \ \thinspace\textbf{MSC Classification: \ }03E70, 03F65} is Constructive
Morse Set Theory (CMST), a formal foundation for Bishop-style constructive
mathematics outlined in \cite{dsb2022} and presented in detail in
\cite{BAlps}. In order that the reader may understand the notation, etc., used
in our paper, and until \cite{BAlps} appears, we refer him/her to
\cite{dsb2022}, \cite{Morse2}, or \cite{Alps}.

CMST\ is based on intuitionistic logic; is a constructive counterpart to A.P.
Morse's classical-logic-based\ \emph{A Theory of Sets} \cite{Morse} (see also
\cite{Morse2,Alps}); uses an updated version of Morse's theory of language and
notation; has a universal set $\mathsf{U}$, membership of which appears to
capture the informal idea of a set being well constructed; and has the
distinctive feature that every mathematical/logical object may be regarded as
either a set or a proposition. This last feature suggested to the author that
we should be able to include (representations of) proofs---or
\emph{justifications,} as we prefer to call them in order to avoid confusion
with the metamathematical notion of proof---as terms of CMST, under additional
axioms such as ones capturing the BHK interpretation\footnote{%
%TCIMACRO{\TeXButton{sf}{\normalfont\sf}}%
%BeginExpansion
\normalfont\sf
%EndExpansion
Brouwer-HeytingKolmogorov interpretation} of intuitionistic logic. This would
provide a form of internal realizability, in comparison with standard forms of
metamathematical realizability (for example, as found in \cite{Beeson},
\cite{KV}, and \cite[Volume I, pages 195--203]{TvD}).%

%TCIMACRO{\TeXButton{medskip}{\medskip}}%
%BeginExpansion
\medskip
%EndExpansion%
%TCIMACRO{\TeXButton{noindent}{\noindent}}%
%BeginExpansion
\noindent
%EndExpansion
In the following, we present an axiomatic approach to such a justification
theory within CMST, together with some consequences of the
justification-theoretic axioms. In particular, we show that the axioms of
intuitionistic logic in \cite{BAlps,dsb2022} (see the Appendix) are validated
in our theory. Our main interest lies in two new notions:
`$P~\mathsf{justifies\ }p$' and `$\mathsf{The\ justification\_set\ of\ }p$'.
Following Morse, before we state our new axioms we present an orienting
definition, as well as the definitions and definitional axioms needed in
addition to those in \cite{BAlps,dsb2022}.%

%TCIMACRO{\TeXButton{medskip}{\medskip}}%
%BeginExpansion
\medskip
%EndExpansion

\subsection{4.0 \ \ Orienting definition}

\begin{enumerate}
\item $((P\ \mathsf{justifies}\ p)\equiv(P\ \mathsf{justifies\ }p))$
\end{enumerate}

\subsection{4.1 \ \ Definitions}

\begin{enumerate}
\item $(\mathsf{jst\,}Pp\equiv(P\ \mathsf{justifies}\ p))$

\item $(\mathsf{J}p\equiv\mathsf{E}X\,\mathsf{jst\,}Xp)$

\item $(\mathsf{The\ justification\_set\ of\ }p\equiv\mathsf{J}p)$
\end{enumerate}

\subsection{4.2 \ \ Definitional axioms for justification theory}

We emphasise that these definitional axioms are in addition to the ones needed
for set theory and described in \cite{BAlps, dsb2022}.

\begin{enumerate}
\item $(\mathsf{The\ }x\ $\underline{$\mathsf{u}$}$x\equiv(\mathsf{One\ }%
x\ $\underline{$\mathsf{u}$}$x\rightarrow%
%TCIMACRO{\tbigwedge }%
%BeginExpansion
{\textstyle\bigwedge}
%EndExpansion
x($\underline{$\mathsf{u}$}$x\rightarrow x)))$

\item $(\mathsf{crd\ }xT\equiv%
%TCIMACRO{\tbigvee }%
%BeginExpansion
{\textstyle\bigvee}
%EndExpansion
y(((x,y)\in T)\wedge y))$

\item $(\mathsf{on\ }A\ \mathsf{to\ }B\ \mathsf{is\ }f\equiv
(\mathsf{function\ is\ }f\wedge(\mathsf{dmn\ }f=A)\wedge(\mathsf{rng\ }%
f\subset B)))$

\item $(.fx\equiv\mathsf{The\ }y((x,y)\in f)))$

\item $(\mathsf{jst\,}Pp\equiv(P~\mathsf{justifies\ }p))$

\item $(\mathsf{J}p\equiv\mathsf{E}X\,\mathsf{jst\,}Xp)$

\item $(\mathsf{dd}\ A\equiv\mathsf{E}x,y((x\in A)\wedge\mathsf{jst\,}y(x\in
A)))$
\end{enumerate}

\subsection{4.3 \ \ Axioms for justification theory}

\begin{enumerate}
\item $(p\leftrightarrow%
%TCIMACRO{\tbigvee }%
%BeginExpansion
{\textstyle\bigvee}
%EndExpansion
P\ \mathsf{jst\,}Pp)$

\item $(\,\mathsf{jst\,}Pp\rightarrow(P\in\mathsf{U}))$

\item $(\mathsf{J}p\in\mathsf{U})$

\item $((\mathsf{jst\ }Pp\wedge\ \mathsf{jst\ }Pq)\rightarrow p\leftrightarrow
q)$

\item $(\mathsf{jst\,}P(p\wedge q)\leftrightarrow\mathsf{orderedpair\ is}%
\ P\wedge\mathsf{jst\,crd}^{\prime}Pp\wedge\mathsf{jst\,crd}^{\prime\prime
}Pq)$

\item $(\mathsf{jst\,}P(p\vee q)\leftrightarrow\mathsf{orderedpair\ is\ }P$

$\ \ \ \ \ \ \ \ \ \ \ \ \ \ \ \ \ \ \ \ \wedge(((\mathsf{crd}^{\prime
}P=0)\wedge\mathsf{jst\,crd}^{\prime\prime}Pp)\vee((\mathsf{crd}^{\prime
}P=1)\wedge\mathsf{jst\,crd}^{\prime\prime}Pq)))$

\item $(\mathsf{jst\,}F(p\rightarrow q)\leftrightarrow\mathsf{on\ J}%
p\ \mathsf{to\ J}q\ \mathsf{is\ }F)$

\item $((A\in\mathsf{U)}\rightarrow\mathsf{jst\,}F%
%TCIMACRO{\tbigwedge }%
%BeginExpansion
{\textstyle\bigwedge}
%EndExpansion
x\in A\,\underline{\mathsf{u}}x\leftrightarrow\mathsf{on\ dd\,}%
A\ \mathsf{is\ }F\wedge%
%TCIMACRO{\tbigwedge }%
%BeginExpansion
{\textstyle\bigwedge}
%EndExpansion
t\in\mathsf{dd\,}A\ \mathsf{jst\,}.Ft\underline{\mathsf{u}}\mathsf{crd}%
^{\prime}t)$

\item $((A\in\mathsf{U)}\rightarrow\mathsf{jst\,}T%
%TCIMACRO{\tbigvee }%
%BeginExpansion
{\textstyle\bigvee}
%EndExpansion
x\in A\,\underline{\mathsf{u}}x\leftrightarrow\mathsf{orderedtriple\ is\ }%
T\wedge\mathsf{jst\,crd}^{\prime\prime}T(\mathsf{crd}^{\prime}T\in A)$

$\ \ \ \ \ \ \ $\ $\ \ $\ $\ \ $\ \ $\ $\ $\ \ \wedge\mathsf{jst\,crd}%
^{\prime\prime\prime}T\underline{\mathsf{u}}\mathsf{crd}^{\prime}T)$

\item $((x\in\mathsf{J}p)\rightarrow\mathsf{One\ }P\ \mathsf{jst\ }%
P(x\in\mathsf{J}p\mathsf{))}$
\end{enumerate}

%

%TCIMACRO{\TeXButton{smallskip}{\smallskip}}%
%BeginExpansion
\smallskip
%EndExpansion%
%TCIMACRO{\TeXButton{noindent}{\noindent}}%
%BeginExpansion
\noindent
%EndExpansion
Axiom (4.3.1)\footnote{%
%TCIMACRO{\TeXButton{sf}{\normalfont\sf}}%
%BeginExpansion
\normalfont\sf
%EndExpansion
Originally our theory was to be Chapter 4 of \cite{BAlps}. It---and the
intended Chapter 3---were omittted on the grounds of the book's size. However,
it is convenient to keep the numberings 4.X for the sections, other than the
introduction, of our paper, and (4.X.Y) for theorem Y of Section X. A
numbering of the form ($\mathsf{A}$.X.Y), with$\ \mathsf{A}=1,2,\mathsf{or\ }%
3$, refers to theorem Y of Section X in Chapter $\mathsf{A}$ of \cite{BAlps}.}
needs no explanation. Axiom (4.3.2) tells us that justifications are
well-constructed terms. Axiom (4.3.4) says that if $P$ justifies two
statements, then those statements are equivalent. Axioms (4.3.5)--(4.3.9)
reflect the BHK interpretation of connectives and quantifiers. We shall return
later to the reason for the hypothesis `$(A\in\mathsf{U})$' in (4.3.8) and (4.3.9).%

%TCIMACRO{\TeXButton{medskip}{\medskip}}%
%BeginExpansion
\medskip
%EndExpansion%
%TCIMACRO{\TeXButton{noindent}{\noindent}}%
%BeginExpansion
\noindent
%EndExpansion
Axiom (4.3.10) may appear strange at first sight. It says that if $x$ is a
justification of $p$, then there is a unique justification showing that $x$ is
a justification of $P$. The underlying informal idea is that the unique
justification in question is simply a check that $x$ is a derivation of $p$
according to the rules of logic. Taken in conjunction with our other axioms
for justification theory, axiom (4.3.10) leads \emph{inter alia} to a
principle of choice from justification sets ((4.4.39) and (4.4.40) below).

\subsection{4.4\ Theorems}

\begin{enumerate}
\item $(P\in\mathsf{J}p\leftrightarrow\mathsf{jst\,}Pp)$%

%TCIMACRO{\TeXButton{smallskip}{\smallskip}}%
%BeginExpansion
\smallskip
%EndExpansion
\textsc{proof}. \ \ By (4.1.2), (2.13.7), and (4.3.2),%
\begin{equation}
(P\in\mathsf{J}p\leftrightarrow P\in\mathsf{E}P\,\mathsf{jst\,}%
Pp\leftrightarrow P\in\mathsf{U}\wedge\mathsf{jst\,}Pp\leftrightarrow
\mathsf{jst\,}Pp). \tag*{$\square$}%
\end{equation}

\item $(\mathsf{J}p\in\mathsf{U})$%

%TCIMACRO{\TeXButton{smallskip}{\smallskip}}%
%BeginExpansion
\smallskip
%EndExpansion
\textsc{proof}.\footnote{%
%TCIMACRO{\TeXButton{sf}{\normalfont\sf}}%
%BeginExpansion
\normalfont\sf
%EndExpansion
This proof is an application of \emph{Neveln's trick}, based on the logical
theorem $(p\leftrightarrow(\mathsf{U}\rightarrow p))$: to prove `$p$', it is
enough to prove `$(\mathsf{U}\rightarrow p)$' \cite[(1.5.6)]{BAlps}.} \ \ By
logic, (4.3.1), (4.3.7), and (2.21.56),%
\begin{align}
(\mathsf{U}  &  \rightarrow(p\rightarrow p)\nonumber\\
&  \rightarrow%
%TCIMACRO{\tbigvee }%
%BeginExpansion
{\textstyle\bigvee}
%EndExpansion
F(F\in\mathsf{J}(p\rightarrow p))\nonumber\\
&  \rightarrow%
%TCIMACRO{\tbigvee }%
%BeginExpansion
{\textstyle\bigvee}
%EndExpansion
F(F\in\mathsf{U}\wedge\mathsf{on\ J}p\ \mathsf{to\ J}p\ \mathsf{is\ }%
f\nonumber\\
&  \rightarrow%
%TCIMACRO{\tbigvee }%
%BeginExpansion
{\textstyle\bigvee}
%EndExpansion
F(F\in\mathsf{U}\wedge\mathsf{function\ is\ }F\wedge\mathsf{dmn\ }%
F=\mathsf{J}p)\nonumber\\
&  \rightarrow%
%TCIMACRO{\tbigvee }%
%BeginExpansion
{\textstyle\bigvee}
%EndExpansion
F(\mathsf{function\ is\ }F\wedge\mathsf{J}p=\mathsf{dmn\ }F\in\mathsf{U}%
)\nonumber\\
&  \rightarrow\mathsf{J}p\in\mathsf{U}). \tag*{$\square$}%
\end{align}

\item $(\lnot p\leftrightarrow\mathsf{J}p=0)$%

%TCIMACRO{\TeXButton{smallskip}{\smallskip}}%
%BeginExpansion
\smallskip
%EndExpansion
\textsc{proof}. \ \ By (4.3.1),%
\begin{equation}
(\lnot p\leftrightarrow\lnot%
%TCIMACRO{\tbigvee }%
%BeginExpansion
{\textstyle\bigvee}
%EndExpansion
x(x\in\mathsf{J}p)\leftrightarrow\mathsf{J}p=0). \tag*{$\square$}%
\end{equation}

\item $(\mathsf{J}0=0)$

\item $(p\rightarrow\mathsf{J}(\lnot p)=0)$%

%TCIMACRO{\TeXButton{smallskip}{\smallskip}}%
%BeginExpansion
\smallskip
%EndExpansion
\textsc{proof}. \ \ By (4.3.1) and (2.6.32),%
\begin{align}
(p  &  \rightarrow(%
%TCIMACRO{\tbigvee }%
%BeginExpansion
{\textstyle\bigvee}
%EndExpansion
x(x\in\mathsf{J(\lnot}p))\rightarrow\lnot p\rightarrow\lnot p\wedge
p\rightarrow0)\nonumber\\
&  \rightarrow\lnot%
%TCIMACRO{\tbigvee }%
%BeginExpansion
{\textstyle\bigvee}
%EndExpansion
x(x\in\mathsf{J(\lnot}p))\nonumber\\
&  \rightarrow\mathsf{J}(\lnot p)=0). \tag*{$\square$}%
\end{align}

\item $(f\in\mathsf{J}(\lnot p)\leftrightarrow\lnot p\wedge f=0)$%

%TCIMACRO{\TeXButton{smallskip}{\smallskip}}%
%BeginExpansion
\smallskip
%EndExpansion
\textsc{proof}. \ \ By (2.9.22), (4.3.7), (4.4.4), (2.21.92), and (4.4.3),%
\begin{align}
(\mathsf{U}  &  \rightarrow\lnot p=(p\rightarrow0)\nonumber\\
&  \rightarrow f\in\mathsf{J}(\lnot p)\leftrightarrow f\in\mathsf{J}%
(p\rightarrow0)\nonumber\\
&  \ \ \ \ \ \ \ \ \ \ \ \ \ \ \ \ \ \,\left.  \leftrightarrow\mathsf{on\ }%
\mathsf{J}p\ \mathsf{to\ }\mathsf{J}0\text{ }\mathsf{is\ }f\wedge
\mathsf{J}0=0\right. \nonumber\\
&  \ \ \ \ \ \ \ \ \ \ \ \ \ \ \ \ \ \,\left.  \leftrightarrow\mathsf{on\ }%
\mathsf{J}p\ \mathsf{to\ }0\text{ }\mathsf{is\ }f\right. \nonumber\\
&  \ \ \ \ \ \ \ \ \ \ \ \ \ \ \ \ \ \,\left.  \leftrightarrow\mathsf{J}%
p=0\wedge f=0\right. \nonumber\\
&  \ \ \ \ \ \ \ \ \ \ \ \ \ \ \ \ \ \,\left.  \leftrightarrow\lnot p\wedge
f=0)\right.  . \tag*{$\square$}%
\end{align}

\item $(\lnot p\leftrightarrow\mathsf{jst\,}0(\lnot p)\leftrightarrow
\mathsf{J}(\lnot p)=1)$%

%TCIMACRO{\TeXButton{smallskip}{\smallskip}}%
%BeginExpansion
\smallskip
%EndExpansion
\textsc{proof}. \ \ By (4.4.6) and (2.11.24),%
\[
(\lnot p\rightarrow%
%TCIMACRO{\tbigwedge }%
%BeginExpansion
{\textstyle\bigwedge}
%EndExpansion
f(f\in\mathsf{J}(\lnot p)\leftrightarrow f=0)\rightarrow\mathsf{J}(\lnot
p)=\{0\}=1)
\]
On the other hand, by (4.3.1),%
\begin{equation}
(\mathsf{J}(\lnot p)=1\rightarrow0\in\mathsf{J}(\lnot p)\rightarrow\lnot p).
\tag*{$\square$}%
\end{equation}

\item $(\mathsf{J}(\lnot p)=0\rightarrow\lnot\lnot p)$%

%TCIMACRO{\TeXButton{smallskip}{\smallskip}}%
%BeginExpansion
\smallskip
%EndExpansion
\textsc{proof}. \ \ By (4.4.7),%
\begin{equation}
(\mathsf{J}(\lnot p)=0\rightarrow\lnot\mathsf{J}(\lnot p)=1\rightarrow
\lnot\lnot p). \tag*{$\square$}%
\end{equation}

\item $(p\rightarrow\mathsf{J}(\lnot\lnot p)=1)$%

%TCIMACRO{\TeXButton{smallskip}{\smallskip}}%
%BeginExpansion
\smallskip
%EndExpansion
\textsc{proof}.\ \ \ By (4.4.5) and (4.4.7),%
\begin{equation}
(p\rightarrow\mathsf{J}(\lnot p)=0\rightarrow\lnot\lnot p\rightarrow
\mathsf{J}(\lnot\lnot p)=1). \tag*{$\square$}%
\end{equation}

\item $(f\in\mathsf{J}(p\rightarrow\lnot\lnot p)\leftrightarrow
f=\rightthreetimes x\in\mathsf{J}p0)$%

%TCIMACRO{\TeXButton{smallskip}{\smallskip}}%
%BeginExpansion
\smallskip
%EndExpansion
\textsc{proof}.\ \ By (4.3.7), (4.3.1), (4.3.10), and (2.20.168),%
\begin{align*}
(f\in\mathsf{J}(p\rightarrow\lnot\lnot p)  &  \rightarrow\mathsf{on\ J}%
p\ \mathsf{to\ J}(\lnot\lnot p)\ \mathsf{is\ }f\\
&  \rightarrow\mathsf{on\ J}p\ \mathsf{is\ }f\wedge%
%TCIMACRO{\tbigwedge }%
%BeginExpansion
{\textstyle\bigwedge}
%EndExpansion
x(x\in\mathsf{J}p\rightarrow.fx\in\mathsf{J}(\lnot\lnot p))\\
&  \rightarrow\mathsf{on\ J}p\ \mathsf{is\ }f\wedge\\
&  \ \ \ \ \left.  \wedge%
%TCIMACRO{\tbigwedge }%
%BeginExpansion
{\textstyle\bigwedge}
%EndExpansion
x(x\in\mathsf{J}p\rightarrow p\rightarrow\mathsf{J}(\lnot\lnot p)=1\rightarrow
.fx=0)\right. \\
&  \rightarrow f=\mathsf{E}x,y(x\in\mathsf{J}p\wedge y=0)=\rightthreetimes
x\in\mathsf{J}p0).
\end{align*}
Conversely, by (2.20.186), (4.4.9), and (4.3.7),%
\begin{align}
(f=\rightthreetimes x\in\mathsf{J}p0  &  \rightarrow\mathsf{function\ is\ }%
f\wedge\mathsf{dmn\ }f=\mathsf{E}x\in\mathsf{J}p(0\in\mathsf{U})=\mathsf{J}%
p\nonumber\\
&  \ \ \ \ \ \left.  \wedge%
%TCIMACRO{\tbigwedge }%
%BeginExpansion
{\textstyle\bigwedge}
%EndExpansion
x\in\mathsf{dmn\ }f(.fx=0)\right. \nonumber\\
&  \rightarrow\mathsf{on\ J}p\ \mathsf{is\ }f\wedge%
%TCIMACRO{\tbigwedge }%
%BeginExpansion
{\textstyle\bigwedge}
%EndExpansion
x(x\in\mathsf{J}p\rightarrow p\rightarrow.fx=0\in1=\mathsf{J}(\lnot\lnot
p))\nonumber\\
&  \rightarrow\mathsf{on\ J}p\ \mathsf{is\ }f\wedge\mathsf{rng\ }%
f\subset\mathsf{J}(\lnot\lnot p)\nonumber\\
&  \rightarrow\mathsf{on\ J}p\ \mathsf{to\ J}(\lnot\lnot p)\mathsf{\ is\ }%
f\nonumber\\
&  \rightarrow f\in\mathsf{J}(p\rightarrow\lnot\lnot p)). \tag*{$\square$}%
\end{align}

\item $(\mathsf{jst\,}(x,y)(p\wedge q)\leftrightarrow\mathsf{jst\,}%
xp\wedge\mathsf{jst\,}yq)$

\item $(\mathsf{jst\,}(m,x)(p\vee q)\leftrightarrow m=0\wedge\mathsf{jst\,}%
xp\vee m=1\wedge\mathsf{jst\,}xq)$

\item $(\mathsf{J}(p\rightarrow q)=\mathsf{maps\ J}p\mathsf{J}q=\mathsf{On\ J}%
p\cap\mathsf{To\ J}q)$

\item $\mathsf{relation\ is\ dd}A$

\item $(x,P\in\mathsf{dd\ }A\leftrightarrow x\in A\wedge\mathsf{jst\ }P(x\in
A))$%

%TCIMACRO{\TeXButton{smallskip}{\smallskip}}%
%BeginExpansion
\smallskip
%EndExpansion
\textsc{proof}. \ \ By (4.3.2), (2.16.77), and (4.2.6),%
\begin{align}
(x\in A\wedge\mathsf{jst\ }P(x\in A)  &  \rightarrow x\in\mathsf{U}\wedge
P\in\mathsf{U}\wedge x\in A\wedge\mathsf{jst\ }P(x\in A)\nonumber\\
&  \rightarrow x,P\in\mathsf{U}\wedge x\in A\wedge\mathsf{jst\ }P(x\in
A)\nonumber\\
&  \rightarrow x,P\in\mathsf{dd}~A). \tag*{$\square$}%
\end{align}

\item $(A\in\mathsf{U\leftrightarrow dd\,}A\in\mathsf{U})$%

%TCIMACRO{\TeXButton{smallskip}{\smallskip}}%
%BeginExpansion
\smallskip
%EndExpansion
\textsc{proof}. \ \ By logic, (4.3.2), (4.3.8), and (2.21.56),%
\begin{align*}
(A\in\mathsf{U}  &  \rightarrow%
%TCIMACRO{\tbigwedge }%
%BeginExpansion
{\textstyle\bigwedge}
%EndExpansion
x\in A(x\in A)\\
&  \rightarrow%
%TCIMACRO{\tbigvee }%
%BeginExpansion
{\textstyle\bigvee}
%EndExpansion
F\ \mathsf{jst\ }F%
%TCIMACRO{\tbigwedge }%
%BeginExpansion
{\textstyle\bigwedge}
%EndExpansion
x\in A(x\in A)\\
&  \rightarrow%
%TCIMACRO{\tbigvee }%
%BeginExpansion
{\textstyle\bigvee}
%EndExpansion
F(F\in\mathsf{U}\wedge\mathsf{jst\,}F%
%TCIMACRO{\tbigwedge }%
%BeginExpansion
{\textstyle\bigwedge}
%EndExpansion
x\in A(x\in A))\\
&  \rightarrow%
%TCIMACRO{\tbigvee }%
%BeginExpansion
{\textstyle\bigvee}
%EndExpansion
F(F\in\mathsf{U}\wedge\mathsf{function\ is\ }F\wedge\mathsf{dmn\ }%
F=\mathsf{dd\,}A)\\
&  \rightarrow%
%TCIMACRO{\tbigvee }%
%BeginExpansion
{\textstyle\bigvee}
%EndExpansion
F(\mathsf{dd\,}A\in\mathsf{U})\\
&  \rightarrow\mathsf{dd\,}A\in\mathsf{U}).
\end{align*}
Conversely, by (4.3.1), (4.4.15), and (2.19.60),%
\begin{align*}
(%
%TCIMACRO{\tbigwedge }%
%BeginExpansion
{\textstyle\bigwedge}
%EndExpansion
x(x\in A  &  \leftrightarrow%
%TCIMACRO{\tbigvee }%
%BeginExpansion
{\textstyle\bigvee}
%EndExpansion
P\,\mathsf{jst\,}P(x\in A)\\
&  \leftrightarrow%
%TCIMACRO{\tbigvee }%
%BeginExpansion
{\textstyle\bigvee}
%EndExpansion
P(x\in A\wedge\mathsf{jst\,}P(x\in A))\\
&  \leftrightarrow%
%TCIMACRO{\tbigvee }%
%BeginExpansion
{\textstyle\bigvee}
%EndExpansion
P(x,P\in\mathsf{dd\,}A)\\
&  \leftrightarrow x\in\mathsf{dmn\,dd\,}A),
\end{align*}
so, by (4.4.14), (4.4.16), and (2.19.77),%
\begin{equation}
(\mathsf{U}\rightarrow\mathsf{relation\ is\ dd}A\wedge A=\mathsf{dmn\,dd\,}%
A\wedge\mathsf{dd}A\in\mathsf{U}\rightarrow A\in\mathsf{U}). \tag*{$\square$}%
\end{equation}

\item $(A\in\mathsf{U}\rightarrow\mathsf{jst\,}F%
%TCIMACRO{\tbigwedge }%
%BeginExpansion
{\textstyle\bigwedge}
%EndExpansion
x\in A\,\underline{\mathsf{u}}x$

$\ \ \ \ \ \ \ \ \ \ \ \ \ \ \ \ \ \leftrightarrow F\in\mathsf{On\ dd\,}%
A\wedge%
%TCIMACRO{\tbigwedge }%
%BeginExpansion
{\textstyle\bigwedge}
%EndExpansion
x,Y(x\in A\wedge\mathsf{jst\,}Y(x\in A)\rightarrow\mathsf{jst\,}%
.F(x,y)\underline{\mathsf{u}}x))$

\item $(A\in\mathsf{U}\rightarrow\mathsf{jst\,}F%
%TCIMACRO{\tbigwedge }%
%BeginExpansion
{\textstyle\bigwedge}
%EndExpansion
x\in A\,\underline{\mathsf{u}}x$

$\ \ \ \ \ \ \ \ \ \ \ \ \ \ \leftrightarrow F\in\mathsf{On\ dd\,}A\wedge%
%TCIMACRO{\tbigwedge }%
%BeginExpansion
{\textstyle\bigwedge}
%EndExpansion
x,Y(\mathsf{jst\ }Y(x\in A)\rightarrow\mathsf{jst\,}%
.F(x,y)\underline{\mathsf{u}}x))$%

%TCIMACRO{\TeXButton{smallskip}{\smallskip}}%
%BeginExpansion
\smallskip
%EndExpansion
\textsc{proof}. \ \ By (4.4.15), (4.4.4), and (2.17.74),%
\begin{align*}
&  \left.  (F\in\mathsf{On\ dd\,}A\right. \\
&  \ \ \ \left.  \rightarrow%
%TCIMACRO{\tbigwedge }%
%BeginExpansion
{\textstyle\bigwedge}
%EndExpansion
x,Y(x,Y\in\mathsf{dd\ }A\leftrightarrow x\in A\wedge\mathsf{jst\ }Y(x\in
A)\leftrightarrow\mathsf{jst\ }Y(x\in A))\right. \\
&  \ \ \ \left.  \rightarrow%
%TCIMACRO{\tbigwedge }%
%BeginExpansion
{\textstyle\bigwedge}
%EndExpansion
x,Y(x,Y\in\mathsf{dd\ }A\rightarrow\mathsf{jst\,}.F(x,y)\underline{\mathsf{u}%
}x)\right. \\
&
\ \ \ \ \ \ \ \ \ \ \ \ \ \ \ \ \ \ \ \ \ \ \ \ \ \ \ \ \ \ \ \ \ \ \ \ \ \left.
\leftrightarrow%
%TCIMACRO{\tbigwedge }%
%BeginExpansion
{\textstyle\bigwedge}
%EndExpansion
x,Y(\mathsf{jst\ }Y(x\in A)\rightarrow\mathsf{jst\,}%
.F(x,y)\underline{\mathsf{u}}x))).\right.
\end{align*}
The result now follows from (4.4.17).%
%TCIMACRO{\TeXButton{hfill}{\hfill}}%
%BeginExpansion
\hfill
%EndExpansion
$\square$

\item $(A\in\mathsf{U}\rightarrow\mathsf{jst\ }F%
%TCIMACRO{\tbigwedge }%
%BeginExpansion
{\textstyle\bigwedge}
%EndExpansion
x\in A\ y\leftrightarrow F\in\mathsf{On\ dd\ }A\ \cap\mathsf{To\ J}y)$%

%TCIMACRO{\TeXButton{smallskip}{\smallskip}}%
%BeginExpansion
\smallskip
%EndExpansion
\textsc{proof}.\ \ \ Taking $(\underline{\mathsf{u}}x=y)$, we see from (4.3.8)
and (4.4.1) that%
\begin{align}
(A\in\mathsf{U}\rightarrow(\mathsf{jst\ }F%
%TCIMACRO{\tbigwedge }%
%BeginExpansion
{\textstyle\bigwedge}
%EndExpansion
x\in A\ y  &  \leftrightarrow F\in\mathsf{On\ dd\ }A\wedge%
%TCIMACRO{\tbigwedge }%
%BeginExpansion
{\textstyle\bigwedge}
%EndExpansion
t\in\mathsf{dd}\ A\ \mathsf{jst\ }.Ft\,\underline{\mathsf{u}}\mathsf{crd}%
^{\prime}t\nonumber\\
&  \leftrightarrow F\in\mathsf{On\ dd\ }A\wedge%
%TCIMACRO{\tbigwedge }%
%BeginExpansion
{\textstyle\bigwedge}
%EndExpansion
t\in\mathsf{dd}\ A\ \mathsf{jst\ }.Ft\,y\nonumber\\
&  \leftrightarrow F\in\mathsf{On\ dd\ }A\wedge%
%TCIMACRO{\tbigwedge }%
%BeginExpansion
{\textstyle\bigwedge}
%EndExpansion
t\in\mathsf{dd}\ A(.Ft\in\mathsf{J}y)\nonumber\\
&  \leftrightarrow F\in\mathsf{On\ dd\ }A\ \cap\mathsf{To\ J}y\ ). \tag*{$
\square$}%
\end{align}

\item $(A\in\mathsf{U}\rightarrow\mathsf{jst\,}T%
%TCIMACRO{\tbigvee }%
%BeginExpansion
{\textstyle\bigvee}
%EndExpansion
x\in A\,y$

$\ \ \ \ \ \ \ \ \ \ \ \ \ \ \ \ \ \ \ \ \ \ \ \ \leftrightarrow
\mathsf{orderedtriple\ is\ }T\wedge\mathsf{jst\,crd}^{\prime\prime
}T(\mathsf{crd}^{\prime}t\in A)\wedge\mathsf{jst\,crd}^{\prime\prime\prime
}Ty)$

\item \label{21}$(A,B,\in\mathsf{U}\rightarrow\mathsf{jst\,}F%
%TCIMACRO{\tbigwedge }%
%BeginExpansion
{\textstyle\bigwedge}
%EndExpansion
x\in A%
%TCIMACRO{\tbigvee }%
%BeginExpansion
{\textstyle\bigvee}
%EndExpansion
y\in B\underline{\mathsf{u}}^{\prime}xy$

$\ \ \ \ \ \ \ \ \ \ \ \ \ \ \ \ \leftrightarrow F\in\mathsf{On\ dd\,}%
A$\ $\wedge\
%TCIMACRO{\tbigwedge }%
%BeginExpansion
{\textstyle\bigwedge}
%EndExpansion
t\in\mathsf{dd\,}A(\mathsf{orderedtriple\ is\ }.Ft$

$\ \ \ \ \ \ \ \ \ \ \ \ \ \ \ \ \ \ \ \ \ \ \ \ \ \ \ \ \ \wedge
\mathsf{jst\ crd}^{\prime\prime}.Ft(\mathsf{crd}^{\prime}.Ft\in B)\wedge
\mathsf{jst\,crd}^{\prime\prime\prime}.Ft\,\underline{\mathsf{u}}^{\prime
}\mathsf{crd}^{\prime}t\,\mathsf{crd}^{\prime}.Ft))$%

%TCIMACRO{\TeXButton{smallskip}{\smallskip}}%
%BeginExpansion
\smallskip
%EndExpansion
\textsc{proof}. \ \ By (4.3.8) and (4.3.9),
\begin{align}
&  (A,B,\in\mathsf{U\rightarrow}\nonumber\\
&  \ \ \left.  \mathsf{jst\,}F%
%TCIMACRO{\tbigwedge }%
%BeginExpansion
{\textstyle\bigwedge}
%EndExpansion
x\in A%
%TCIMACRO{\tbigvee }%
%BeginExpansion
{\textstyle\bigvee}
%EndExpansion
y\in B\ \underline{\mathsf{u}}^{\prime}xy\right. \nonumber\\
&  \ \ \ \ \left.  \leftrightarrow F\in\mathsf{On\ dd\ }A\wedge%
%TCIMACRO{\tbigwedge }%
%BeginExpansion
{\textstyle\bigwedge}
%EndExpansion
t\in\mathsf{dd\,}A\ \mathsf{jst\,}.Ft%
%TCIMACRO{\tbigvee }%
%BeginExpansion
{\textstyle\bigvee}
%EndExpansion
y\in B\underline{\mathsf{u}}^{\prime}\mathsf{crd\,}^{\prime}ty\right.
\nonumber\\
&  \ \ \ \ \left.  \leftrightarrow F\in\mathsf{On\ dd\ }A\right. \nonumber\\
&  \ \ \ \ \ \ \ \ \left.  \wedge\
%TCIMACRO{\tbigwedge }%
%BeginExpansion
{\textstyle\bigwedge}
%EndExpansion
t\in\mathsf{dd\,}A(\mathsf{orderedtriple\ is\ }.Ft\wedge\mathsf{jst\ crd}%
^{\prime\prime}.Ft(\mathsf{crd}^{\prime}.Ft\in B)\right. \nonumber\\
&  \ \ \ \ \ \ \ \ \ \ \ \ \ \ \ \ \ \ \ \ \ \ \left.  \wedge
\ \mathsf{jst\,crd}^{\prime\prime\prime}.Ft\,\underline{\mathsf{u}}^{\prime
}\mathsf{crd}^{\prime}t\,\mathsf{crd}^{\prime}.Ft)\right.  \tag{$\square$}%
\end{align}

\end{enumerate}

%

%TCIMACRO{\TeXButton{medskip}{\medskip}}%
%BeginExpansion
\medskip
%EndExpansion%
%TCIMACRO{\TeXButton{noindent}{\noindent}}%
%BeginExpansion
\noindent
%EndExpansion
We now translate the axioms of logic into the language of our justification theory.

\begin{enumerate}
\item[22.] $(F=\rightthreetimes Q\in\mathsf{J}q\rightthreetimes P\in
\mathsf{J}p\ Q\rightarrow\mathsf{jst\ }F(q\rightarrow(p\rightarrow q)))$%

%TCIMACRO{\TeXButton{smallskip}{\smallskip}}%
%BeginExpansion
\smallskip
%EndExpansion
\textsc{proof}.\ \ By (2.20.186) and (4.3.7),%
\begin{align*}
&  \left.
%TCIMACRO{\tbigwedge }%
%BeginExpansion
{\textstyle\bigwedge}
%EndExpansion
Q(Q\in\mathsf{J}q\wedge g=\rightthreetimes P\in\mathsf{J}p\ Q\right. \\
&  \ \ \ \left.  \rightarrow\mathsf{function\ is\ }g\wedge\mathsf{dmn\ }%
g=\mathsf{E}P\in\mathsf{J}p(Q\in\mathsf{U})=\mathsf{J}p\right. \\
&  \ \ \ \ \ \ \ \ \left.  \wedge\,%
%TCIMACRO{\tbigwedge }%
%BeginExpansion
{\textstyle\bigwedge}
%EndExpansion
P\in\mathsf{J}p(.gP=Q)\right. \\
&  \ \ \ \left.  \rightarrow\mathsf{on\ J}p\ \mathsf{to\ J}q\text{
}\mathsf{is}\text{ }g\right. \\
&  \ \ \ \left.  \rightarrow g\in\mathsf{J}(p\rightarrow q)).\right.
\end{align*}
Hence, again by (2.20.186) and (4.3.7),%
\begin{align}
&  \left.  (F=\rightthreetimes Q\in\mathsf{J}q\rightthreetimes P\in
\mathsf{J}p\ Q\right. \nonumber\\
&  \ \ \ \ \ \left.  \rightarrow\mathsf{function\ is\ }F\wedge\mathsf{dmn\ }%
F=\mathsf{E}Q\in\mathsf{J}q((\rightthreetimes P\in\mathsf{J}p\,Q)\in
\mathsf{U})\right. \nonumber\\
&  \ \ \ \ \ \ \ \ \ \ \left.  \wedge\,%
%TCIMACRO{\tbigwedge }%
%BeginExpansion
{\textstyle\bigwedge}
%EndExpansion
Q\in\mathsf{J}q(.FQ=(\rightthreetimes P\in\mathsf{J}p\,Q)\in\mathsf{J}%
(p\rightarrow q))\right. \nonumber\\
&  \ \ \ \ \ \left.  \rightarrow\mathsf{on\ J}q\ \mathsf{to\ J}(p\rightarrow
q)\mathsf{\ is\ }F\right. \nonumber\\
&  \ \ \ \ \ \left.  \rightarrow\mathsf{jst\ }F(q\rightarrow(p\rightarrow
q))).\right.  \tag*{$\square$}%
\end{align}

\item[23.] $(\mathsf{jst\,}X(p\rightarrow(q\rightarrow r))$

\ \ \ $\rightarrow(F=\rightthreetimes Y\in\mathsf{J}(p\rightarrow
q)\rightthreetimes P\in\mathsf{J}p(..XP.YP)$

$\ \ \ \ \ \ \ \ \ \ \ \ \ \ \rightarrow\mathsf{jst\,}F((p\rightarrow
q)\rightarrow(p\rightarrow r))))$%

%TCIMACRO{\TeXButton{smallskip}{\smallskip}}%
%BeginExpansion
\smallskip
%EndExpansion
\textsc{proof}.\ \ \ Let%
\begin{align*}
(\alpha &  =(\mathsf{jst\,}X(p\rightarrow(q\rightarrow r))\wedge
\mathsf{jst\,}Y(p\rightarrow q))\ \wedge\\
F  &  =\rightthreetimes Y\in\mathsf{J}(p\rightarrow q)\rightthreetimes
P\in\mathsf{J}p(..XP.YP)).
\end{align*}
By (4.3.7),%
\begin{align*}
(\alpha\rightarrow%
%TCIMACRO{\tbigwedge }%
%BeginExpansion
{\textstyle\bigwedge}
%EndExpansion
P(P\in\mathsf{J}p  &  \rightarrow\mathsf{on\ J}p\ \mathsf{to\ }\mathsf{J}%
(q\rightarrow r)\mathsf{\ is\ }X\wedge\mathsf{on\ }\mathsf{J}p\ \mathsf{to\ }%
\mathsf{J}q\ \mathsf{is\ }Y\\
&  \rightarrow.XP\in\mathsf{J}(q\rightarrow r)\wedge.YP\in\mathsf{J}q\\
&  \rightarrow\mathsf{on\ J}q\ \mathsf{to\ J}r\ \mathsf{is\ }.XP\wedge
.YP\in\mathsf{J}q\\
&  \rightarrow..XP.YP\in\mathsf{J}r)).
\end{align*}
Hence, by (2.20.186), (4.3.7), and (4.4.1),%
\begin{align*}
&  \ \ \left.  \ (\alpha\wedge g=\rightthreetimes P\in\mathsf{J}%
p(..XP.YP)\right. \\
&  \ \ \ \ \ \left.  \rightarrow\mathsf{function\ is\ }g\wedge\mathsf{dmn\ }%
g=\mathsf{E}P\in\mathsf{J}p(..XP.YP\in\mathsf{U})=\mathsf{J}p\right. \\
&  \ \ \ \ \ \ \ \ \left.  \wedge\
%TCIMACRO{\tbigwedge }%
%BeginExpansion
{\textstyle\bigwedge}
%EndExpansion
P\in\mathsf{dmn\ }g(.gP=..XP.YP\in\mathsf{J}r)\right. \\
&  \ \ \ \ \ \left.  \rightarrow\mathsf{on\ J}p\ \mathsf{to\ J}%
r\ \mathsf{is\ }g\right. \\
&  \ \ \ \ \ \left.  \rightarrow\mathsf{jst\ }g(p\rightarrow r)\right. \\
&  \ \ \ \ \ \left.  \rightarrow g\in\mathsf{J}(p\rightarrow r)).\right.
\end{align*}
Again using (2.20.186) and (4.3.7), we have%
\begin{align}
&  \left.  (\mathsf{jst\,}X(p\rightarrow(q\rightarrow r))\right. \nonumber\\
&  \ \ \ \ \left.  \rightarrow\mathsf{function\ is\ }F\right. \nonumber\\
&  \ \ \ \ \ \ \ \ \ \left.  \wedge\ \mathsf{dmn\ }F=\mathsf{E}Y\in
\mathsf{J}(p\rightarrow q)((\rightthreetimes P\in\mathsf{J}p(..XP.YP))\in
\mathsf{U})=\mathsf{J}(p\rightarrow q)\right. \nonumber\\
&  \ \ \ \ \ \ \ \ \ \left.  \wedge%
%TCIMACRO{\tbigwedge }%
%BeginExpansion
{\textstyle\bigwedge}
%EndExpansion
Y\in\mathsf{J}(p\rightarrow q)(.FY=(\rightthreetimes P\in\mathsf{J}%
p(..XP.YP))\in\mathsf{J}(p\rightarrow r))\right. \nonumber\\
&  \ \ \ \ \left.  \rightarrow\mathsf{on\ J}(p\rightarrow q)\ \mathsf{to\ J}%
(p\rightarrow r)~\mathsf{is\ }F\right. \nonumber\\
&  \ \ \ \ \left.  \rightarrow\mathsf{jst\ }F((p\rightarrow q)\rightarrow
(p\rightarrow r))).\right.  \tag*{$\square$}%
\end{align}

\item[24.] $(\mathsf{jst\,}X(p\wedge q)\rightarrow\mathsf{jst\,crd}^{\prime
}X\,p)$

\item[25.] $(\mathsf{jst\,}X(p\wedge q)\rightarrow\mathsf{jst\,crd}%
^{\prime\prime}X\,q)$

\item[26.] $(\mathsf{jst\,}Xp\rightarrow(F=\rightthreetimes Y\in
\mathsf{J}q(X,Y)\rightarrow\mathsf{jst\ }F(q\rightarrow p\wedge q))$%

%TCIMACRO{\TeXButton{smallskip}{\smallskip}}%
%BeginExpansion
\smallskip
%EndExpansion
\textsc{proof}. \ \ By (4.3.5), (4.4.2), (2.20.186), and (4.3.7),%
\begin{align}
(\mathsf{jst\ }Xp  &  \rightarrow%
%TCIMACRO{\tbigwedge }%
%BeginExpansion
{\textstyle\bigwedge}
%EndExpansion
Y(Y\in\mathsf{J}q\rightarrow\mathsf{orderedpair\ is\ }X,Y\wedge\mathsf{jst\,}%
\mathsf{crd}^{\prime}(X,Y)p\wedge\mathsf{jst\,}\mathsf{crd}^{\prime\prime
}(X,Y)q\nonumber\\
&  \ \ \ \ \ \ \ \ \ \ \ \ \ \ \ \ \ \ \ \,\left.  \rightarrow\mathsf{jst\,}%
(X,Y)(p\wedge q)\right. \nonumber\\
&  \ \ \ \ \ \ \ \ \ \ \ \ \ \ \ \ \ \ \ \,\left.  \rightarrow X,Y\in
\mathsf{J}(p\wedge q)))\right. \nonumber\\
&  \rightarrow(F=\rightthreetimes Y\in\mathsf{J}q(X,Y)\rightarrow
\mathsf{function\ is\ }F\wedge\mathsf{dmn\ }F=\mathsf{J}q\nonumber\\
&
\ \ \ \ \ \ \ \ \ \ \ \ \ \ \ \ \ \ \ \ \ \ \ \ \ \ \ \ \ \ \ \ \ \ \ \ \left.
\wedge%
%TCIMACRO{\tbigwedge }%
%BeginExpansion
{\textstyle\bigwedge}
%EndExpansion
Y\in\mathsf{J}q(.FY=X,Y\in\mathsf{J}(p\wedge q))\right. \nonumber\\
&  \rightarrow\mathsf{on\ J}q\ \mathsf{to\ J}(p\wedge q)~\mathsf{is\ }%
F\nonumber\\
&  \rightarrow\mathsf{jst\ }F(q\rightarrow p\wedge q)). \tag*{$\square$}%
\end{align}

\item[27.] $(\mathsf{jst\,}Xp\rightarrow\mathsf{jst\,}(0,X)(p\vee q))$

\item[28.] $(\mathsf{jst\,}Xq\rightarrow\mathsf{jst\,}(1,X)(p\vee q))$

\item[29.] $(\mathsf{jst\,}X((p\rightarrow r)\wedge(q\rightarrow r))$

$\ \ \ \rightarrow$\ $(F=\rightthreetimes Y\in\mathsf{J}(p\vee q)(\mathsf{crd}%
^{\prime}Y=0\wedge.\mathsf{crd}^{\prime}X\mathsf{crd}^{\prime\prime}Y$

$\ \ \ \ \ \ \ \ \ \ \ \ \ \ \ \ \ \ \ \ \ \ \ \ \ \ \ \ \ \ \ \ \ \!\vee
\mathsf{crd}^{\prime}Y=1\wedge.\mathsf{crd}^{\prime\prime}X\mathsf{crd}%
^{\prime\prime}Y)\rightarrow\mathsf{jst\,}F(p\vee q\rightarrow r)))$%

%TCIMACRO{\TeXButton{smallskip}{\smallskip}}%
%BeginExpansion
\smallskip
%EndExpansion
\textsc{proof}. \ \ By (4.3.5), (4.3.6), and (4.3.7),%
\begin{align*}
&  \left.  (\mathsf{jst\,}X((p\rightarrow r)\wedge(q\rightarrow r))\wedge
Y\in\mathsf{J}(p\vee q)\right. \\
&  \ \ \ \ \left.  \rightarrow\mathsf{orderedpair\ is\ }X\wedge\mathsf{jst\,}%
\mathsf{crd}^{\prime}X(p\rightarrow r)\wedge\mathsf{jst\,}\mathsf{crd}%
^{\prime\prime}X(q\rightarrow r)\right. \\
&  \ \ \ \ \ \ \ \ \left.  \wedge\ \mathsf{orderedpair\ is\ }Y\right. \\
&  \ \ \ \ \ \ \ \ \left.  \wedge\ (\mathsf{crd}^{\prime}Y=0\wedge
\mathsf{jst\,}\mathsf{crd}^{\prime\prime}Yp\vee\mathsf{crd}^{\prime}%
Y=1\wedge\mathsf{jst\,}\mathsf{crd}^{\prime\prime}Yq)\right. \\
&  \ \ \ \ \left.  \rightarrow\mathsf{on\ }\mathsf{J}p\ \mathsf{to\ }%
\mathsf{J}r\text{ }\mathsf{is\ crd}^{\prime}X\wedge\mathsf{on\ }%
\mathsf{J}q\ \mathsf{to\ }\mathsf{J}r\text{ }\mathsf{is\ crd}^{\prime\prime
}X\right. \\
&  \ \ \ \ \ \ \ \ \left.  \wedge\ (\mathsf{crd}^{\prime}Y=0\wedge
\mathsf{crd}^{\prime\prime}Y\in\mathsf{J}p\vee\mathsf{crd}^{\prime}%
Y=1\wedge\mathsf{crd}^{\prime\prime}Y\in\mathsf{J}q)\right. \\
&  \ \ \ \ \left.  \rightarrow\mathsf{on\ }\mathsf{J}p\ \mathsf{to\ }%
\mathsf{J}r\text{ }\mathsf{is\ crd}^{\prime}X\wedge\mathsf{on\ }%
\mathsf{J}q\ \mathsf{to\ }\mathsf{J}r\text{ }\mathsf{is\ crd}^{\prime\prime
}X\right. \\
&  \ \ \ \ \ \ \ \ \left.  \wedge\ (\mathsf{crd}^{\prime}Y=0\wedge
.\mathsf{crd}^{\prime}X\mathsf{crd}^{\prime\prime}Y\in\mathsf{J}%
r\vee\mathsf{crd}^{\prime}Y=1\wedge.\mathsf{crd}^{\prime\prime}X\mathsf{crd}%
^{\prime\prime}Y\in\mathsf{J}r))\right.
\end{align*}
Thus%
\begin{align*}
&  \left.  (\mathsf{jst\,}X((p\rightarrow r)\wedge(q\rightarrow r))\right. \\
&  \left.  \rightarrow%
%TCIMACRO{\tbigwedge }%
%BeginExpansion
{\textstyle\bigwedge}
%EndExpansion
Y\in\mathsf{J}(p\vee q)(\mathsf{crd}^{\prime}Y=0\wedge.\mathsf{crd}^{\prime
}X\mathsf{crd}^{\prime\prime}Y\in\mathsf{J}r\right. \\
&  \ \ \ \ \ \ \ \ \ \ \ \ \ \ \ \ \ \ \ \ \left.  \vee\ \,\mathsf{crd}%
^{\prime}Y=1\wedge.\mathsf{crd}^{\prime\prime}X\mathsf{crd}^{\prime\prime}%
Y\in\mathsf{J}r))\right.
\end{align*}
and therefore, by (2.20.186) and (4.3.7),%
\begin{align}
&  \left.  (\mathsf{jst\,}X((p\rightarrow r)\wedge(q\rightarrow r))\right.
\nonumber\\
&  \left.  \rightarrow(F=\rightthreetimes Y\in\mathsf{J}(p\vee q)(\mathsf{crd}%
^{\prime}Y=0\wedge.\mathsf{crd}^{\prime}X\mathsf{crd}^{\prime\prime}Y\right.
\nonumber\\
&  \ \ \ \ \ \ \ \ \ \ \ \ \ \ \ \ \ \ \ \ \ \ \ \ \ \ \ \ \left.
\vee\,\,\mathsf{crd}^{\prime}Y=1\wedge.\mathsf{crd}^{\prime\prime
}X\mathsf{crd}^{\prime\prime}Y)\right. \nonumber\\
&  \ \ \ \ \ \ \ \ \left.  \rightarrow\mathsf{on\ J}(p\vee q)\ \mathsf{to\ }%
\mathsf{J}r\ \mathsf{is\ }F\right. \nonumber\\
&  \ \ \ \ \ \ \ \ \left.  \rightarrow\mathsf{jst\,}F(p\vee q\rightarrow
r)))\right.  \tag*{$\square$}%
\end{align}

\end{enumerate}

%

%TCIMACRO{\TeXButton{noindent}{\noindent}}%
%BeginExpansion
\noindent
%EndExpansion
And now for the predicate axioms.\footnote{%
%TCIMACRO{\TeXButton{sf}{\normalfont\sf}}%
%BeginExpansion
\normalfont\sf
%EndExpansion
The hypothesis `$(A\in\mathsf{U})$' in (4.10.33), (4.10.35), (4.10.37), and
(4.10.38) reflects the same hypothesis in (4.3.7), and does not appear in the
actual predicate axioms (which are, after all, logical, not set-theoretic) in
\cite{BAlps, dsb2022}.}%

%TCIMACRO{\TeXButton{smallskip}{\smallskip}}%
%BeginExpansion
\smallskip
%EndExpansion

\begin{enumerate}
\item[30.] $(A\in\mathsf{U\wedge jst\,}P(a\in A)\wedge X=\rightthreetimes
F\in\mathsf{J}%
%TCIMACRO{\tbigwedge }%
%BeginExpansion
{\textstyle\bigwedge}
%EndExpansion
x\in A\,\underline{\mathsf{u}}x\,.F(a,P)$

$\ \ \ \rightarrow\mathsf{j\mathsf{st\,}}X(%
%TCIMACRO{\tbigwedge }%
%BeginExpansion
{\textstyle\bigwedge}
%EndExpansion
x\in A\,\underline{\mathsf{u}}x\rightarrow\mathsf{\underline{\mathsf{u}}}a)).$%

%TCIMACRO{\TeXButton{smallskip}{\smallskip}}%
%BeginExpansion
\smallskip
%EndExpansion
\textsc{proof}. \ \ By (4.3.1), (4.4.15), (4.3.8), and (4.4.1),%
\begin{align*}
&  \left.  (A\in\mathsf{U}\wedge\mathsf{jst\,}P(a\in A)\wedge\mathsf{jst\,}F%
%TCIMACRO{\tbigwedge }%
%BeginExpansion
{\textstyle\bigwedge}
%EndExpansion
x\in A\,\underline{\mathsf{u}}x\right. \\
&  \ \ \ \ \ \left.  \rightarrow a\in A\wedge\mathsf{jst\,}P(a\in A)\right. \\
&  \ \ \ \ \ \left.  \rightarrow a,P\in\mathsf{dd\ }A\wedge\mathsf{on\ dd\,}%
A\ \mathsf{is\ }F\wedge%
%TCIMACRO{\tbigwedge }%
%BeginExpansion
{\textstyle\bigwedge}
%EndExpansion
t\in\mathsf{dd\,}A\ \mathsf{jst\,}.Ft\underline{\mathsf{u}}\mathsf{crd}%
^{\prime}t\right. \\
&  \ \ \ \ \ \left.  \rightarrow\mathsf{jst\ }.F(a,P)\underline{\mathsf{u}%
}a)\right. \\
&  \ \ \ \ \ \left.  \rightarrow.F(a,P)\in\mathsf{J}\underline{\mathsf{u}%
}a)\right.
\end{align*}
From this, (2.20.186), and (4.3.7) we obtain%
\begin{align}
&  \left.  (A\in\mathsf{U\wedge jst\,}P(a\in A)\wedge X=\rightthreetimes
F\in\mathsf{J}%
%TCIMACRO{\tbigwedge }%
%BeginExpansion
{\textstyle\bigwedge}
%EndExpansion
x\in A\,\underline{\mathsf{u}}x\,.F(a,P)\right. \nonumber\\
&  \ \ \ \left.  \rightarrow%
%TCIMACRO{\tbigwedge }%
%BeginExpansion
{\textstyle\bigwedge}
%EndExpansion
F(F\in\mathsf{J}%
%TCIMACRO{\tbigwedge }%
%BeginExpansion
{\textstyle\bigwedge}
%EndExpansion
x\in A\,\underline{\mathsf{u}}x\rightarrow.F(a,P)\in\mathsf{J}%
\underline{\mathsf{u}}a)\right. \nonumber\\
&  \ \ \ \left.  \rightarrow\mathsf{function\ is\ }X\wedge\mathsf{dmn\ }%
X=\mathsf{J}%
%TCIMACRO{\tbigwedge }%
%BeginExpansion
{\textstyle\bigwedge}
%EndExpansion
x\in A\,\underline{\mathsf{u}}x\right. \nonumber\\
&  \ \ \ \ \ \ \ \left.  \wedge%
%TCIMACRO{\tbigwedge }%
%BeginExpansion
{\textstyle\bigwedge}
%EndExpansion
F\in\mathsf{J}%
%TCIMACRO{\tbigwedge }%
%BeginExpansion
{\textstyle\bigwedge}
%EndExpansion
x\in A\,\underline{\mathsf{u}}x(.XF=.F(a,P)\in\mathsf{J}\underline{\mathsf{u}%
}a)\right. \nonumber\\
&  \ \ \ \left.  \rightarrow\mathsf{on\ J}%
%TCIMACRO{\tbigwedge }%
%BeginExpansion
{\textstyle\bigwedge}
%EndExpansion
x\in A\,\underline{\mathsf{u}}x\ \mathsf{to\ J}\underline{\mathsf{u}%
}a\ \mathsf{is\ }X\right. \nonumber\\
&  \ \ \ \left.  \rightarrow\mathsf{j\mathsf{st\,}}X(%
%TCIMACRO{\tbigwedge }%
%BeginExpansion
{\textstyle\bigwedge}
%EndExpansion
x\in A\,\underline{\mathsf{u}}x\rightarrow\mathsf{\underline{\mathsf{u}}%
}a)).\right.  \tag*{$\square$}%
\end{align}

\item[31.] $(A\in\mathsf{U}\wedge\mathsf{jst\,}P(a\in A)\wedge
F=\rightthreetimes Q\in\mathsf{J}\underline{\mathsf{u}}a\,(a,P,Q)\rightarrow
\mathsf{jst\,}F(\underline{\mathsf{u}}a\rightarrow%
%TCIMACRO{\tbigvee }%
%BeginExpansion
{\textstyle\bigvee}
%EndExpansion
x\in A\,\underline{\mathsf{u}}x)))$%

%TCIMACRO{\TeXButton{smallskip}{\smallskip}}%
%BeginExpansion
\smallskip
%EndExpansion
\textsc{proof}. \ \ By (2.16.80), (4.3.2), (2.16.74), and (4.3.9),
\begin{align*}
&  \left.  (A\in\mathsf{U}\wedge\mathsf{jst\,}P(a\in A)\wedge Q\in
\mathsf{J}\underline{\mathsf{u}}a\wedge T=(a,P,Q)\right. \\
&  \ \ \ \left.  \rightarrow A\in\mathsf{U}\wedge a\in A\wedge P\in
\mathsf{U}\wedge Q\in\mathsf{U}\right. \\
&  \ \ \ \left.  \rightarrow T\in\mathsf{U}\wedge\mathsf{orderedtriple\ is\ }%
T\wedge\ \mathsf{jst\,}P(a\in A)\right. \\
&  \ \ \ \ \ \ \ \ \left.  \wedge~\mathsf{jst\ }Q\underline{\mathsf{u}}%
a\wedge\ \mathsf{crd}^{\prime}T=a\wedge\ \mathsf{crd}^{\prime\prime}%
T=P\wedge\ \mathsf{crd}^{\prime\prime}T=Q\right. \\
&  \ \ \ \left.  \rightarrow T\in\mathsf{U}\wedge\mathsf{orderedtriple\ is\ }%
T\wedge\ \mathsf{jst\ crd}^{\prime\prime}T(\mathsf{crd}^{\prime}T\in
A)\wedge\ \mathsf{jst\ crd}^{\prime\prime\prime}T\underline{\mathsf{u}%
}\mathsf{crd}^{\prime}T\right.  \ \\
&  \ \ \ \left.  \rightarrow\mathsf{jst\ }T%
%TCIMACRO{\tbigvee }%
%BeginExpansion
{\textstyle\bigvee}
%EndExpansion
x\in A\underline{\mathsf{u}}x\right. \\
&  \ \ \ \left.  \rightarrow T\in\mathsf{J}%
%TCIMACRO{\tbigvee }%
%BeginExpansion
{\textstyle\bigvee}
%EndExpansion
x\in A\underline{\mathsf{u}}x).\right.
\end{align*}
Hence, by (2.20.186) and (4.3.7),%
\begin{align}
&  \left.  (A\in\mathsf{U}\wedge\mathsf{jst\,}P(a\in A)\wedge
F=\rightthreetimes Q\in\mathsf{J}\underline{\mathsf{u}}a\,(a,P,Q)\right.
\nonumber\\
&  \ \ \ \left.  \rightarrow%
%TCIMACRO{\tbigwedge }%
%BeginExpansion
{\textstyle\bigwedge}
%EndExpansion
Q(Q\in\mathsf{J}\underline{\mathsf{u}}a\rightarrow a,P,Q\in\mathsf{J}%
%TCIMACRO{\tbigvee }%
%BeginExpansion
{\textstyle\bigvee}
%EndExpansion
x\in A\underline{\mathsf{u}}x)\right. \nonumber\\
&  \ \ \ \left.  \rightarrow\mathsf{function\ is\ }F\wedge\mathsf{dmn\ }%
F=\mathsf{E}x\in\mathsf{J}\underline{\mathsf{u}}a(a,P,Q\in\mathsf{U})\right.
\nonumber\\
&  \ \ \ \ \ \ \ \left.  \wedge%
%TCIMACRO{\tbigwedge }%
%BeginExpansion
{\textstyle\bigwedge}
%EndExpansion
Q\in\mathsf{J}\underline{\mathsf{u}}a(.FQ=a,P,Q\in\mathsf{J}%
%TCIMACRO{\tbigvee }%
%BeginExpansion
{\textstyle\bigvee}
%EndExpansion
x\in A\underline{\mathsf{u}}x)\right. \nonumber\\
&  \ \ \ \left.  \rightarrow\mathsf{on}\ \mathsf{J}\underline{\mathsf{u}%
}a\ \mathsf{to\ J}%
%TCIMACRO{\tbigvee }%
%BeginExpansion
{\textstyle\bigvee}
%EndExpansion
x\in A\underline{\mathsf{u}}x\ \mathsf{is\ }F\right. \nonumber\\
&  \ \ \ \left.  \rightarrow\mathsf{jst\ }F(\underline{\mathsf{u}}a\rightarrow%
%TCIMACRO{\tbigvee }%
%BeginExpansion
{\textstyle\bigvee}
%EndExpansion
x\in A\,\underline{\mathsf{u}}x)).\right.  \tag*{$\square$}%
\end{align}

\item[32.] $(A\in\mathsf{U}\wedge F=\rightthreetimes P\in\mathsf{J}%
y\rightthreetimes t\in\mathsf{dd\,}A\,P\rightarrow\mathsf{jst\,}F(y\rightarrow%
%TCIMACRO{\tbigwedge }%
%BeginExpansion
{\textstyle\bigwedge}
%EndExpansion
x\in A\,y))$%

%TCIMACRO{\TeXButton{smallskip}{\smallskip}}%
%BeginExpansion
\smallskip
%EndExpansion
\textsc{proof}. \ \ By (2.20.186), (4.3.2), (4.3.8), and (4.4.1),%
\begin{align*}
&  \left.  (A\in\mathsf{U}\wedge P\in\mathsf{J}y\wedge f=\rightthreetimes
t\in\mathsf{dd\,}A\,P\right. \\
&  \ \ \ \left.  \rightarrow\mathsf{function\ is\ }f\wedge\mathsf{dmn\ }%
f=\mathsf{E}t\in\mathsf{dd}\ A(P\in\mathsf{U})=\mathsf{dd\ }A\right. \\
&  \ \ \ \ \ \ \ \left.  \wedge\
%TCIMACRO{\tbigwedge }%
%BeginExpansion
{\textstyle\bigwedge}
%EndExpansion
t\in\mathsf{dmn\ }f(.ft=P)\right. \\
&  \ \ \ \left.  \rightarrow\mathsf{on\ dd\ }A\ \mathsf{to\ J}y\ \mathsf{is\ }%
f\wedge%
%TCIMACRO{\tbigwedge }%
%BeginExpansion
{\textstyle\bigwedge}
%EndExpansion
t\in\mathsf{dd}A\,\mathsf{jst\ }.fty\right. \\
&  \ \ \ \left.  \rightarrow\mathsf{jst\ }f%
%TCIMACRO{\tbigwedge }%
%BeginExpansion
{\textstyle\bigwedge}
%EndExpansion
x\in A\,y\right. \\
&  \ \ \ \left.  \rightarrow f\in\mathsf{J}%
%TCIMACRO{\tbigwedge }%
%BeginExpansion
{\textstyle\bigwedge}
%EndExpansion
x\in A\,y).\right.
\end{align*}
Hence, by (2.20.186) and (4.3.7),%
\begin{align}
&  \left.  (A\in\mathsf{U}\wedge F=\rightthreetimes P\in\mathsf{J}%
y\rightthreetimes t\in\mathsf{dd\,}A\,P\right. \nonumber\\
&  \ \ \ \left.  \rightarrow\mathsf{function\ is\ }F\wedge\mathsf{dmn\ }%
F=\mathsf{E}P\in\mathsf{J}y(\rightthreetimes t\in\mathsf{dd\,}A\,P\in
\mathsf{U})=\mathsf{J}y\right. \nonumber\\
&  \ \ \ \ \ \ \ \left.  \wedge%
%TCIMACRO{\tbigwedge }%
%BeginExpansion
{\textstyle\bigwedge}
%EndExpansion
P\in\mathsf{dmn\ }F(.FP=\rightthreetimes t\in\mathsf{dd\,}A\,P\in\mathsf{J}%
%TCIMACRO{\tbigwedge }%
%BeginExpansion
{\textstyle\bigwedge}
%EndExpansion
x\in A\,y)\right. \nonumber\\
&  \ \ \ \left.  \rightarrow\mathsf{on\ J}y\ \mathsf{to\ J}%
%TCIMACRO{\tbigwedge }%
%BeginExpansion
{\textstyle\bigwedge}
%EndExpansion
x\in A\,y\ \mathsf{is}\ F\right. \nonumber\\
&  \ \ \ \left.  \rightarrow\mathsf{jst\ }F(y\rightarrow%
%TCIMACRO{\tbigwedge }%
%BeginExpansion
{\textstyle\bigwedge}
%EndExpansion
x\in A\,y)).\right.  \tag*{$\square$}%
\end{align}

\item[33.] $(A\in\mathsf{U}\wedge F=\rightthreetimes T\in\mathsf{J}%
%TCIMACRO{\tbigvee }%
%BeginExpansion
{\textstyle\bigvee}
%EndExpansion
x\in A\,y\ \mathsf{crd}^{\prime\prime}T\rightarrow\mathsf{jst\,}F(%
%TCIMACRO{\tbigvee }%
%BeginExpansion
{\textstyle\bigvee}
%EndExpansion
x\in A\,y\rightarrow y))$%

%TCIMACRO{\TeXButton{smallskip}{\smallskip}}%
%BeginExpansion
\smallskip
%EndExpansion
\textsc{proof}.\ \ \ By (4.4.1), (4.3.9), (2.20.186), and (4.3.7),%
\begin{align}
&  \left.  (A\in\mathsf{U}\right. \nonumber\\
&  \ \ \ \left.  \rightarrow%
%TCIMACRO{\tbigwedge }%
%BeginExpansion
{\textstyle\bigwedge}
%EndExpansion
T(T\in\mathsf{J}%
%TCIMACRO{\tbigvee }%
%BeginExpansion
{\textstyle\bigvee}
%EndExpansion
x\in A\,y\right. \nonumber\\
&  \ \ \ \ \ \ \ \ \ \ \ \ \ \left.  \rightarrow\mathsf{jst\ }T%
%TCIMACRO{\tbigvee }%
%BeginExpansion
{\textstyle\bigvee}
%EndExpansion
x\in A\,y\right. \nonumber\\
&  \ \ \ \ \ \ \ \ \ \ \ \ \ \left.  \rightarrow\mathsf{orderedtriple\ is\ }%
T\wedge\mathsf{jst\ crd}^{\prime\prime}T(\mathsf{crd}^{\prime}T\in
A)\wedge\mathsf{jst\ crd}^{\prime\prime}Ty)\right. \nonumber\\
&  \ \ \ \left.  \rightarrow(F=\rightthreetimes T\in\mathsf{J}%
%TCIMACRO{\tbigvee }%
%BeginExpansion
{\textstyle\bigvee}
%EndExpansion
x\in A\,y\ \mathsf{crd}^{\prime\prime}T\right. \nonumber\\
&  \ \ \ \ \ \ \ \ \ \ \left.  \rightarrow\mathsf{function\ is\ }F\right.
\nonumber\\
&  \ \ \ \ \ \ \ \ \ \ \ \ \ \ \ \ \left.  \wedge\ \mathsf{dmn\ }%
F=\mathsf{E}T\in\mathsf{J}%
%TCIMACRO{\tbigvee }%
%BeginExpansion
{\textstyle\bigvee}
%EndExpansion
x\in A\,y(\mathsf{crd}^{\prime\prime}T\in\mathsf{U})=\mathsf{J}%
%TCIMACRO{\tbigvee }%
%BeginExpansion
{\textstyle\bigvee}
%EndExpansion
x\in A\,y\right. \nonumber\\
&  \ \ \ \ \ \ \ \ \ \ \ \ \ \ \ \ \left.  \wedge\
%TCIMACRO{\tbigwedge }%
%BeginExpansion
{\textstyle\bigwedge}
%EndExpansion
T\in\mathsf{dmn\ }F(.FT=\mathsf{crd}^{\prime\prime}T\in\mathsf{J}y)\right.
\nonumber\\
&  \ \ \ \ \ \ \ \ \ \ \left.  \rightarrow\mathsf{on\ J}%
%TCIMACRO{\tbigvee }%
%BeginExpansion
{\textstyle\bigvee}
%EndExpansion
x\in A\,y\ \mathsf{to\ J}y\ \mathsf{is\ }F\right. \nonumber\\
&  \ \ \ \ \ \ \ \ \ \ \left.  \rightarrow\mathsf{jst\,}F(%
%TCIMACRO{\tbigvee }%
%BeginExpansion
{\textstyle\bigvee}
%EndExpansion
x\in A\,y\rightarrow y)).\right.  \tag*{$\square$}%
\end{align}

\item[34.] $(A\in\mathsf{U}\wedge\mathsf{jst\,}F%
%TCIMACRO{\tbigwedge }%
%BeginExpansion
{\textstyle\bigwedge}
%EndExpansion
x\in A(\underline{\mathsf{u}}x\rightarrow\underline{\mathsf{v}}x)$

\ $\wedge\ G=\rightthreetimes\Phi\in\mathsf{J}%
%TCIMACRO{\tbigwedge }%
%BeginExpansion
{\textstyle\bigwedge}
%EndExpansion
x\in A\,\underline{\mathsf{u}}x\rightthreetimes t\in\mathsf{dd\,}A..Ft.\Phi t$

$\ \ \ \ \ \ \ \ \ \rightarrow\mathsf{jst\,}G(%
%TCIMACRO{\tbigwedge }%
%BeginExpansion
{\textstyle\bigwedge}
%EndExpansion
x\in A\,\underline{\mathsf{u}}x\rightarrow%
%TCIMACRO{\tbigwedge }%
%BeginExpansion
{\textstyle\bigwedge}
%EndExpansion
x\in A\,\underline{\mathsf{v}}x)))$%

%TCIMACRO{\TeXButton{smallskip}{\smallskip}}%
%BeginExpansion
\smallskip
%EndExpansion
\textsc{proof}. \ \ Let \ \
\begin{align*}
(\alpha &  =(A\in\mathsf{U}\wedge\mathsf{jst\,}F%
%TCIMACRO{\tbigwedge }%
%BeginExpansion
{\textstyle\bigwedge}
%EndExpansion
x\in A(\underline{\mathsf{u}}x\rightarrow\underline{\mathsf{v}}x))\wedge\\
G  &  =\rightthreetimes\Phi\in\mathsf{J}%
%TCIMACRO{\tbigwedge }%
%BeginExpansion
{\textstyle\bigwedge}
%EndExpansion
x\in A\,\underline{\mathsf{u}}x\rightthreetimes t\in\mathsf{dd\,}A..Ft.\Phi
t).
\end{align*}
By (4.4.1), (4.3.8), (4.4.1), (4.3.7), and (2.20.186),%
\begin{align*}
&  \left.  (\alpha\wedge\Phi\in\mathsf{J}%
%TCIMACRO{\tbigwedge }%
%BeginExpansion
{\textstyle\bigwedge}
%EndExpansion
x\in A\,\underline{\mathsf{u}}x\wedge h=\rightthreetimes t\in\mathsf{dd\,}%
A..Ft.\Phi t\right. \\
&  \ \ \ \left.  \rightarrow\mathsf{on\ dd\,}A\ \mathsf{is\ }F\wedge%
%TCIMACRO{\tbigwedge }%
%BeginExpansion
{\textstyle\bigwedge}
%EndExpansion
t\in\mathsf{dd\,}A\ \mathsf{jst\,}.Ft(\underline{\mathsf{u}}\mathsf{crd}%
^{\prime}t\rightarrow\underline{\mathsf{v}}\mathsf{crd}^{\prime}t)\right. \\
&  \ \ \ \ \ \ \ \left.  \wedge\ \mathsf{on\ dd\,}A\ \mathsf{is\ }\Phi\ \wedge%
%TCIMACRO{\tbigwedge }%
%BeginExpansion
{\textstyle\bigwedge}
%EndExpansion
t\in\mathsf{dd\,}A\ \mathsf{jst\ }.\Phi t\underline{\mathsf{u}}\mathsf{crd}%
^{\prime}t\right. \\
&  \ \ \ \left.  \rightarrow\mathsf{on\ dd\,}A\ \mathsf{is\ }F\wedge%
%TCIMACRO{\tbigwedge }%
%BeginExpansion
{\textstyle\bigwedge}
%EndExpansion
t\in\mathsf{dd\,}A(.Ft\in\mathsf{J}(\underline{\mathsf{u}}\mathsf{crd}%
^{\prime}t\rightarrow\underline{\mathsf{v}}\mathsf{crd}^{\prime}t))\right. \\
&  \ \ \ \ \ \ \ \left.  \wedge\ \mathsf{on\ dd\,}A\ \mathsf{is\ }\Phi\wedge%
%TCIMACRO{\tbigwedge }%
%BeginExpansion
{\textstyle\bigwedge}
%EndExpansion
t\in\mathsf{dd\,}A(.\Phi t\in\mathsf{J}\underline{\mathsf{u}}\mathsf{crd}%
^{\prime}t)\right. \\
&  \ \ \ \left.  \rightarrow%
%TCIMACRO{\tbigwedge }%
%BeginExpansion
{\textstyle\bigwedge}
%EndExpansion
t(t\in\mathsf{dd\,}A\rightarrow\mathsf{on\ J}\underline{\mathsf{u}%
}\mathsf{crd}^{\prime}t\ \mathsf{to\ J}\underline{\mathsf{v}}\mathsf{crd}%
^{\prime}t\ \mathsf{is\ }.Ft\wedge.\Phi t\in\mathsf{J}\underline{\mathsf{u}%
}\mathsf{crd}^{\prime}t)\right. \\
&  \ \ \ \ \ \ \ \ \ \ \ \ \ \ \ \ \ \ \ \ \ \ \ \left.  \rightarrow..Ft.\Phi
t\in\mathsf{J}\underline{\mathsf{v}}\mathsf{crd}^{\prime}t)\right. \\
&  \ \ \ \left.  \rightarrow\mathsf{function\ is\ }h\wedge\mathsf{dmn\ }%
h=\mathsf{E}t\in\mathsf{dd\,}A(..Ft.\Phi t\in\mathsf{U})=\mathsf{dd\,}A\right.
\\
&  \ \ \ \ \ \ \ \left.  \wedge%
%TCIMACRO{\tbigwedge }%
%BeginExpansion
{\textstyle\bigwedge}
%EndExpansion
t\in\mathsf{dmn\ }h(.ht=..Ft.\Phi t\in\mathsf{J}\underline{\mathsf{v}%
}\mathsf{crd}^{\prime}t)\right. \\
&  \ \ \ \left.  \rightarrow A\in\mathsf{U}\wedge\mathsf{on\ dd\,}%
A\ \mathsf{is\ }h\wedge%
%TCIMACRO{\tbigwedge }%
%BeginExpansion
{\textstyle\bigwedge}
%EndExpansion
t\in\mathsf{dd\,}A(\mathsf{jst}.ht\underline{\mathsf{v}}\mathsf{crd}^{\prime
}t)\right. \\
&  \ \ \ \left.  \rightarrow\mathsf{jst\ }h%
%TCIMACRO{\tbigwedge }%
%BeginExpansion
{\textstyle\bigwedge}
%EndExpansion
x\in A\,\underline{\mathsf{v}}x).\right. \\
&  \ \ \ \left.  \rightarrow h\in\mathsf{J}%
%TCIMACRO{\tbigwedge }%
%BeginExpansion
{\textstyle\bigwedge}
%EndExpansion
x\in A\,\underline{\mathsf{v}}x).\right.
\end{align*}
Hence, by logic, (2.20.186), and (4.3.7),%
\begin{align}
(\alpha &  \rightarrow%
%TCIMACRO{\tbigwedge }%
%BeginExpansion
{\textstyle\bigwedge}
%EndExpansion
\Phi\in\mathsf{J}%
%TCIMACRO{\tbigwedge }%
%BeginExpansion
{\textstyle\bigwedge}
%EndExpansion
x\in A\,\underline{\mathsf{u}}x(\rightthreetimes t\in\mathsf{dd\,}A\,..Ft.\Phi
t\in\mathsf{J}%
%TCIMACRO{\tbigwedge }%
%BeginExpansion
{\textstyle\bigwedge}
%EndExpansion
x\in A\,\underline{\mathsf{v}}x)\nonumber\\
&  \rightarrow\mathsf{on\ }\mathsf{J}%
%TCIMACRO{\tbigwedge }%
%BeginExpansion
{\textstyle\bigwedge}
%EndExpansion
x\in A\,\underline{\mathsf{u}}x\ \mathsf{to\ }\mathsf{J}%
%TCIMACRO{\tbigwedge }%
%BeginExpansion
{\textstyle\bigwedge}
%EndExpansion
x\in A\,\underline{\mathsf{v}}x\mathsf{\ is\ }G\nonumber\\
&  \rightarrow\mathsf{jst\,}G(%
%TCIMACRO{\tbigwedge }%
%BeginExpansion
{\textstyle\bigwedge}
%EndExpansion
x\in A\,\underline{\mathsf{u}}x\rightarrow%
%TCIMACRO{\tbigwedge }%
%BeginExpansion
{\textstyle\bigwedge}
%EndExpansion
x\in A\,\underline{\mathsf{v}}x)). \tag*{$\square$}%
\end{align}

\item[35.] $(A\in\mathsf{U}\wedge\mathsf{jst\,}F%
%TCIMACRO{\tbigwedge }%
%BeginExpansion
{\textstyle\bigwedge}
%EndExpansion
x\in A(\underline{\mathsf{u}}x\rightarrow\underline{\mathsf{v}}x)$

$\ \ \ \wedge G=\rightthreetimes T\in\mathsf{J}%
%TCIMACRO{\tbigvee }%
%BeginExpansion
{\textstyle\bigvee}
%EndExpansion
x\in A\,\underline{\mathsf{u}}x(\mathsf{crd}^{\prime}T,\mathsf{crd}%
^{\prime\prime}T,.F(\mathsf{crd}^{\prime}T,\mathsf{crd}^{\prime}%
\mathsf{crd}^{\prime\prime}T))$

\ \ \ \ \ \ \ \ \ \ \ \ $\rightarrow\mathsf{jst\,}G(%
%TCIMACRO{\tbigvee }%
%BeginExpansion
{\textstyle\bigvee}
%EndExpansion
x\in A\,\underline{\mathsf{u}}x\rightarrow%
%TCIMACRO{\tbigvee }%
%BeginExpansion
{\textstyle\bigvee}
%EndExpansion
x\in A\,\underline{\mathsf{v}}x))$%

%TCIMACRO{\TeXButton{smallskip}{\smallskip}}%
%BeginExpansion
\smallskip
%EndExpansion
\textsc{proof}. \ \ Let%
\[
(\alpha=(A\in\mathsf{U}\wedge\mathsf{jst\,}F%
%TCIMACRO{\tbigwedge }%
%BeginExpansion
{\textstyle\bigwedge}
%EndExpansion
x\in A(\underline{\mathsf{u}}x\rightarrow\underline{\mathsf{v}}x)))
\]
and%
\[%
%TCIMACRO{\tbigwedge }%
%BeginExpansion
{\textstyle\bigwedge}
%EndExpansion
T(\underline{\mathsf{w}}T=(\mathsf{crd}^{\prime}T,\mathsf{crd}^{\prime\prime
}T,.F(\mathsf{crd}^{\prime}T,\mathsf{crd}^{\prime}\mathsf{crd}^{\prime\prime
}T)).
\]
\label{0102}Using (4.4.1), (4.3.9), (4.3.1), (4.4.15), and (4.3.7), we obtain%
\begin{align*}
&  \left.  (\alpha\wedge T\in\mathsf{J}%
%TCIMACRO{\tbigvee }%
%BeginExpansion
{\textstyle\bigvee}
%EndExpansion
x\in A\,\underline{\mathsf{u}}x\right. \\
&  \ \ \ \left.  \rightarrow\mathsf{jst\ }T%
%TCIMACRO{\tbigvee }%
%BeginExpansion
{\textstyle\bigvee}
%EndExpansion
x\in A\,\underline{\mathsf{u}}x)\right. \\
&  \ \ \ \left.  \rightarrow\mathsf{orderedtriple\ is\ }T\wedge\mathsf{jst\,}%
\mathsf{crd}^{\prime\prime}T(\mathsf{crd}^{\prime}T\in A)\wedge
\mathsf{jst\ crd}^{\prime\prime\prime}T\underline{\mathsf{u}}\mathsf{crd}%
^{\prime}T)\right. \\
&  \ \ \ \left.  \rightarrow\mathsf{orderedtriple\ is\ }T\wedge\mathsf{crd}%
^{\prime}T\in A\wedge\mathsf{jst\,}\mathsf{crd}^{\prime\prime}T(\mathsf{crd}%
^{\prime}T\in A)\wedge\mathsf{crd}^{\prime\prime\prime}T\in\mathsf{J}%
\underline{\mathsf{u}}\mathsf{crd}^{\prime}T)\right. \\
&  \ \ \ \left.  \rightarrow(\mathsf{crd}^{\prime}T,\mathsf{crd}^{\prime
\prime}T)\in\mathsf{dd}A\wedge\underline{\mathsf{u}}\mathsf{crd}^{\prime
}T\wedge\mathsf{on\ dd}A\ \mathsf{is\ }F\right. \\
&  \ \ \ \ \ \ \ \ \left.  \wedge%
%TCIMACRO{\tbigwedge }%
%BeginExpansion
{\textstyle\bigwedge}
%EndExpansion
t\in\mathsf{dd}A\ \mathsf{jst\ }.Ft(\underline{\mathsf{u}}\mathsf{crd}%
^{\prime}t\rightarrow\underline{\mathsf{v}}\mathsf{crd}^{\prime}t)\right. \\
&  \ \ \ \left.  \rightarrow\mathsf{jst}.F(\mathsf{crd}^{\prime}%
T,\mathsf{crd}^{\prime\prime}T)(\underline{\mathsf{u}}\mathsf{crd}^{\prime
}T\rightarrow\underline{\mathsf{v}}\mathsf{crd}^{\prime}T)\wedge\right. \\
&  \ \ \ \left.  \rightarrow\mathsf{on\ J}\underline{\mathsf{u}}%
\mathsf{crd}^{\prime}T\ \mathsf{to\ J}\underline{\mathsf{v}}\mathsf{crd}%
^{\prime}T\ \mathsf{is\ }.F(\mathsf{crd}^{\prime}T,\mathsf{crd}^{\prime\prime
}T)\right. \\
&  \ \ \ \left.  \rightarrow..F(\mathsf{crd}^{\prime}T,\mathsf{crd}%
^{\prime\prime}T)\mathsf{crd}^{\prime\prime\prime}T\in\mathsf{J}%
\underline{\mathsf{v}}\mathsf{crd}^{\prime}T\right. \\
&  \ \ \ \left.  \rightarrow\mathsf{orderedtriple\ is\ }\underline{\mathsf{w}%
}T\wedge\mathsf{jst\ crd}^{\prime\prime}\underline{\mathsf{w}}T(\mathsf{crd}%
^{\prime}\underline{\mathsf{w}}T\in A)\wedge\mathsf{jst\ crd}^{\prime
\prime\prime}\underline{\mathsf{w}}T\underline{\mathsf{v}}\mathsf{crd}%
^{\prime}T\right. \\
&  \ \ \ \left.  \rightarrow\underline{\mathsf{w}}T\in\mathsf{J}%
%TCIMACRO{\tbigvee }%
%BeginExpansion
{\textstyle\bigvee}
%EndExpansion
x\in A\,\underline{\mathsf{v}}x).\right.
\end{align*}
Hence, by logic, (2.20.186), and (4.3.7),%
\begin{align}
&  \left.  (\alpha\wedge G=\rightthreetimes T\in\mathsf{J}%
%TCIMACRO{\tbigvee }%
%BeginExpansion
{\textstyle\bigvee}
%EndExpansion
x\in A\,\underline{\mathsf{u}}x\,\underline{\mathsf{w}}T\right. \nonumber\\
&  \ \ \ \ \ \left.  \rightarrow%
%TCIMACRO{\tbigwedge }%
%BeginExpansion
{\textstyle\bigwedge}
%EndExpansion
T\in\mathsf{J}%
%TCIMACRO{\tbigvee }%
%BeginExpansion
{\textstyle\bigvee}
%EndExpansion
x\in A\,\underline{\mathsf{u}}x(\mathsf{\underline{\mathsf{w}}}T\in\mathsf{J}%
%TCIMACRO{\tbigvee }%
%BeginExpansion
{\textstyle\bigvee}
%EndExpansion
x\in A\ \underline{\mathsf{v}}x)\right. \nonumber\\
&  \ \ \ \ \ \left.  \rightarrow\mathsf{function\ is\ }G\right. \nonumber\\
\ \  &  \ \ \ \ \ \ \ \ \ \left.  \wedge\ \mathsf{dmn\ }G=\mathsf{E}%
T\in\mathsf{J}%
%TCIMACRO{\tbigvee }%
%BeginExpansion
{\textstyle\bigvee}
%EndExpansion
x\in A\,\underline{\mathsf{u}}x(\underline{\mathsf{w}}T\in\mathsf{U)=J}%
%TCIMACRO{\tbigvee }%
%BeginExpansion
{\textstyle\bigvee}
%EndExpansion
x\in A\,\underline{\mathsf{u}}x\right. \nonumber\\
\ \ \  &  \ \ \ \ \ \ \ \ \ \left.  \wedge%
%TCIMACRO{\tbigwedge }%
%BeginExpansion
{\textstyle\bigwedge}
%EndExpansion
T\in\mathsf{J}%
%TCIMACRO{\tbigvee }%
%BeginExpansion
{\textstyle\bigvee}
%EndExpansion
x\in A\,\underline{\mathsf{u}}x(.GT=\underline{\mathsf{w}}T\in\mathsf{J}%
%TCIMACRO{\tbigvee }%
%BeginExpansion
{\textstyle\bigvee}
%EndExpansion
x\in A\ \underline{\mathsf{v}}x)\right. \nonumber\\
&  \ \ \ \ \ \left.  \rightarrow\mathsf{function\ is\ }G\right. \nonumber\\
&  \ \ \ \ \ \ \ \ \ \left.  \wedge%
%TCIMACRO{\tbigwedge }%
%BeginExpansion
{\textstyle\bigwedge}
%EndExpansion
T\in\mathsf{J}%
%TCIMACRO{\tbigvee }%
%BeginExpansion
{\textstyle\bigvee}
%EndExpansion
x\in A\,\underline{\mathsf{u}}x(.GT\in\mathsf{J}%
%TCIMACRO{\tbigvee }%
%BeginExpansion
{\textstyle\bigvee}
%EndExpansion
x\in A\ \underline{\mathsf{v}}x)\right. \nonumber\\
&  \ \ \ \ \ \left.  \rightarrow\mathsf{on\ J}%
%TCIMACRO{\tbigvee }%
%BeginExpansion
{\textstyle\bigvee}
%EndExpansion
x\in A\,\underline{\mathsf{u}}x\ \mathsf{to\ \mathsf{J}}%
%TCIMACRO{\tbigvee }%
%BeginExpansion
{\textstyle\bigvee}
%EndExpansion
x\in A\mathsf{\ \underline{\mathsf{v}}}x\mathsf{\ is\ }G\right. \nonumber\\
&  \ \ \ \ \ \left.  \rightarrow\mathsf{jst\,}G(%
%TCIMACRO{\tbigvee }%
%BeginExpansion
{\textstyle\bigvee}
%EndExpansion
x\in A\,\underline{\mathsf{u}}x\rightarrow%
%TCIMACRO{\tbigvee }%
%BeginExpansion
{\textstyle\bigvee}
%EndExpansion
x\in A\,\underline{\mathsf{v}}x)).\right.  \tag*{$\square$}%
\end{align}

\end{enumerate}

%

%TCIMACRO{\TeXButton{smallskip}{\smallskip}}%
%BeginExpansion
\smallskip
%EndExpansion%
%TCIMACRO{\TeXButton{noindent}{\noindent}}%
%BeginExpansion
\noindent
%EndExpansion
We have not yet used axiom (4.3.10). In order to do so, we need some
preliminary lemmas.

\begin{enumerate}
\item[36.] $(x,P\in\mathsf{dd\,J}p)\leftrightarrow x\in\mathsf{J}p\wedge
P=\mathsf{The\ }P\ \mathsf{jst\,}P(x\in\mathsf{J}p))$

\item[37.] $(\phi=\rightthreetimes x\in\mathsf{J}p\ \mathsf{The\ }%
P\ \mathsf{jst\ }P(x\in\mathsf{J}p)$

\ \ \ $\rightarrow\phi\in\mathsf{On\ J}p\wedge%
%TCIMACRO{\tbigwedge }%
%BeginExpansion
{\textstyle\bigwedge}
%EndExpansion
x\in\mathsf{J}p(.\phi x=\mathsf{The\ }P\ \mathsf{jst\ }P(x\in\mathsf{J}p)))$

\item[38.] $(\phi=\rightthreetimes x\in\mathsf{J}p\ \mathsf{The\ }%
P\ \mathsf{jst\ }P(x\in\mathsf{J}p)\rightarrow x\in\mathsf{J}p\leftrightarrow
x,.\phi x\in\mathsf{dd\,J}p)$
\end{enumerate}

%

%TCIMACRO{\TeXButton{medskip}{\medskip}}%
%BeginExpansion
\medskip
%EndExpansion%
%TCIMACRO{\TeXButton{noindent}{\noindent}}%
%BeginExpansion
\noindent
%EndExpansion
It is well known that the full Axiom of Choice, normally expressed something like%

\[
(%
%TCIMACRO{\tbigwedge }%
%BeginExpansion
{\textstyle\bigwedge}
%EndExpansion
x\in A%
%TCIMACRO{\tbigvee }%
%BeginExpansion
{\textstyle\bigvee}
%EndExpansion
y\in B\ \underline{\mathsf{u}}^{\prime}xy\rightarrow%
%TCIMACRO{\tbigvee }%
%BeginExpansion
{\textstyle\bigvee}
%EndExpansion
f(\mathsf{on}\ A\ \mathsf{to\ }B\ \mathsf{is\ }f\wedge%
%TCIMACRO{\tbigwedge }%
%BeginExpansion
{\textstyle\bigwedge}
%EndExpansion
x\in A\ \underline{\mathsf{u}}^{\prime}x.fx),
\]%
%TCIMACRO{\TeXButton{smallskip}{\smallskip}}%
%BeginExpansion
\smallskip
%EndExpansion
implies the Law of Excluded Middle; see \cite{BAlps,Diac,Goodmy}. Axiom
(4.3.10) enables us to produce two versions of a principle of choice from
justification sets, the first version being explicitly in the setting of our
axiomatic justification theory.

\begin{enumerate}
\item[39.] $(A,B,\in\mathsf{U}\rightarrow\mathsf{jst\,}F%
%TCIMACRO{\tbigwedge }%
%BeginExpansion
{\textstyle\bigwedge}
%EndExpansion
x\in A%
%TCIMACRO{\tbigvee }%
%BeginExpansion
{\textstyle\bigvee}
%EndExpansion
y\in B\underline{\mathsf{u}}^{\prime}xy$

$\ \ \ \ \ \ \ \ \ \ \ \ \ \ \ \ \leftrightarrow F\in\mathsf{On\ dd\,}%
A$\ $\wedge\
%TCIMACRO{\tbigwedge }%
%BeginExpansion
{\textstyle\bigwedge}
%EndExpansion
t\in\mathsf{dd\,}A(\mathsf{orderedtriple\ is\ }.Ft$

$\ \ \ \ \ \ \ \ \ \ \ \ \ \ \ \ \ \ \ \ \ \ \ \ \ \ \ \ \ \wedge
\mathsf{jst\ crd}^{\prime\prime}.Ft(\mathsf{crd}^{\prime}.Ft\in B)\wedge
\mathsf{jst\,crd}^{\prime\prime\prime}.Ft\,\underline{\mathsf{u}}^{\prime
}\mathsf{crd}^{\prime}t\,\mathsf{crd}^{\prime}.Ft))$

\item[40.] \label{def}$(B\in\mathsf{U}\wedge\mathsf{jst\ }F%
%TCIMACRO{\tbigwedge }%
%BeginExpansion
{\textstyle\bigwedge}
%EndExpansion
x\in\mathsf{J}p%
%TCIMACRO{\tbigvee }%
%BeginExpansion
{\textstyle\bigvee}
%EndExpansion
y\in B\,\underline{\mathsf{u}}^{\prime}xy\wedge\phi=\rightthreetimes
x\in\mathsf{J}p\ \mathsf{The\ }P\ \mathsf{jst\ }P(x\in\mathsf{J}p)$

$\ \ \ \rightarrow(f=\rightthreetimes x\in\mathsf{J}p\ \mathsf{crd}%
^{\prime\prime\prime}.F(x,\mathsf{.}\phi x)\rightarrow f\in\mathsf{On\ J}%
p\cap\mathsf{To\ }B\wedge%
%TCIMACRO{\tbigwedge }%
%BeginExpansion
{\textstyle\bigwedge}
%EndExpansion
x\in\mathsf{J}p\ \underline{\mathsf{u}}^{\prime}x.fx))$%

%TCIMACRO{\TeXButton{smallskip}{\smallskip}}%
%BeginExpansion
\smallskip
%EndExpansion
\textsc{proof}.\ \ \ By \label{00ere1702}(4.4.2), (4.4.21), (4.4.38),
(2.20.186), (4.3.3), and (2.20.64),%
\begin{align}
&  \left.  (B\in\mathsf{U}\wedge\mathsf{jst\ }F%
%TCIMACRO{\tbigwedge }%
%BeginExpansion
{\textstyle\bigwedge}
%EndExpansion
x\in\mathsf{J}p%
%TCIMACRO{\tbigvee }%
%BeginExpansion
{\textstyle\bigvee}
%EndExpansion
y\in B\,\underline{\mathsf{u}}^{\prime}xy\wedge\phi=\rightthreetimes
x\in\mathsf{J}p\ \mathsf{The\ }P\ \mathsf{jst\ }P(x\in\mathsf{J}p)\right.
\nonumber\\
&  \ \ \left.  \rightarrow\mathsf{J}p,B,\in\mathsf{U}\right. \nonumber\\
&  \ \ \left.  \rightarrow F\in\mathsf{On\ dd\,J}p\wedge%
%TCIMACRO{\tbigwedge }%
%BeginExpansion
{\textstyle\bigwedge}
%EndExpansion
t\in\mathsf{dd\,J}p(\mathsf{orderedtriple\ is\ }.Ft\right. \nonumber\\
&  \ \ \ \ \ \ \ \left.  \wedge\ \mathsf{jst\,}\mathsf{crd}^{\prime\prime
}.Ft\,(\mathsf{crd}^{\prime}.Ft\in B)\wedge\mathsf{jst\ crd}^{\prime
\prime\prime}.Ft\ \underline{\mathsf{u}}^{\prime}\mathsf{crd}^{\prime
}t\mathsf{crd}^{\prime}.Ft)\right. \nonumber\\
&  \ \ \left.  \rightarrow%
%TCIMACRO{\tbigwedge }%
%BeginExpansion
{\textstyle\bigwedge}
%EndExpansion
x(x\in\mathsf{J}p\rightarrow x,.\phi x\in\mathsf{dd~J}p\right. \nonumber\\
&  \ \ \ \ \ \ \ \ \ \ \ \ \ \ \ \ \ \ \ \ \left.  \rightarrow
\mathsf{orderedtriple\ is\ }.F(x,.\phi x)\wedge\mathsf{crd}^{\prime}.F(x,.\phi
x)\in B\right. \nonumber\\
&  \ \ \ \ \ \ \ \ \ \ \ \ \ \ \ \ \ \ \ \ \ \ \ \ \ \left.  \wedge
\ \mathsf{jst\ crd}^{\prime\prime\prime}.F(x,.\phi x)\,\underline{\mathsf{u}%
}x\mathsf{crd}^{\prime}.F(x,.\phi x)\right. \nonumber\\
&  \ \ \left.  \rightarrow(f=\rightthreetimes x\in\mathsf{J}p\ \mathsf{crd}%
^{\prime}.F(x,\mathsf{.}\phi x)\right. \nonumber\\
&  \ \ \ \ \ \ \ \ \ \ \ \ \ \left.  \rightarrow\mathsf{function\ is\ }%
f\wedge\mathsf{dmn}\ f=\mathsf{E}x\in\mathsf{J}p(\mathsf{crd}^{\prime
}.F(x,\mathsf{.}\phi x)\in\mathsf{U})\right. \nonumber\\
&  \ \ \ \ \ \ \ \ \ \ \ \ \ \ \ \ \ \ \left.  \wedge\
%TCIMACRO{\tbigwedge }%
%BeginExpansion
{\textstyle\bigwedge}
%EndExpansion
x\in\mathsf{dmn\ }f(.fx=\ \mathsf{crd}^{\prime}.F(x,\mathsf{.}\phi x)\in
B)\right. \nonumber\\
&  \ \ \ \ \ \ \ \ \ \ \ \ \ \ \ \ \ \ \left.  \wedge\
%TCIMACRO{\tbigwedge }%
%BeginExpansion
{\textstyle\bigwedge}
%EndExpansion
x(x\in\mathsf{J}p\rightarrow\mathsf{jst\ crd}^{\prime\prime\prime}.F(x,.\phi
x)\,\underline{\mathsf{u}}x\mathsf{crd}^{\prime}.F(x,.\phi x))\right.
\nonumber\\
&  \ \ \ \ \ \ \ \ \ \ \ \ \ \left.  \rightarrow\mathsf{function\ is\ }%
f\wedge\mathsf{dmn}\ f=\mathsf{J}p\in\mathsf{U\wedge rng\ }f=B\right.
\nonumber\\
&  \ \ \ \ \ \ \ \ \ \ \ \ \ \ \ \ \ \ \left.  \wedge%
%TCIMACRO{\tbigwedge }%
%BeginExpansion
{\textstyle\bigwedge}
%EndExpansion
x\in\mathsf{J}p(.fx=\mathsf{crd}^{\prime}.F(x,\mathsf{.}\phi x)\wedge
\underline{\mathsf{u}}^{\prime}x\mathsf{crd}^{\prime}.F(x,\mathsf{.}\phi
x))\right. \nonumber\\
&  \ \ \ \ \ \ \ \ \ \ \ \ \ \left.  \rightarrow f\in\mathsf{U\wedge
on\ J}p\mathsf{\ to\ }B\ \mathsf{is\ }f\wedge%
%TCIMACRO{\tbigwedge }%
%BeginExpansion
{\textstyle\bigwedge}
%EndExpansion
x\in\mathsf{J}p\,\underline{\mathsf{u}}^{\prime}x.fx\right. \nonumber\\
&  \ \ \ \ \ \ \ \ \ \ \ \ \ \left.  \rightarrow f\in\mathsf{On\ J}%
p\cap\mathsf{To\ }B\wedge%
%TCIMACRO{\tbigwedge }%
%BeginExpansion
{\textstyle\bigwedge}
%EndExpansion
x\in\mathsf{J}p\,\underline{\mathsf{u}}^{\prime}x.fx).\right.  \tag*{$
\square$}%
\end{align}

\item[41.] $(B\in\mathsf{U}\wedge%
%TCIMACRO{\tbigwedge }%
%BeginExpansion
{\textstyle\bigwedge}
%EndExpansion
x\in\mathsf{J}p%
%TCIMACRO{\tbigvee }%
%BeginExpansion
{\textstyle\bigvee}
%EndExpansion
y\in B\,\underline{\mathsf{u}}^{\prime}xy\rightarrow%
%TCIMACRO{\tbigvee }%
%BeginExpansion
{\textstyle\bigvee}
%EndExpansion
f\in\mathsf{On\ J}p\cap\mathsf{To\ }B%
%TCIMACRO{\tbigwedge }%
%BeginExpansion
{\textstyle\bigwedge}
%EndExpansion
x\in\mathsf{J}p\,\underline{\mathsf{u}}^{\prime}x.fx))$
\end{enumerate}

%

%TCIMACRO{\TeXButton{medskip}{\medskip}}%
%BeginExpansion
\medskip
%EndExpansion%
%TCIMACRO{\TeXButton{noindent}{\noindent}}%
%BeginExpansion
\noindent
%EndExpansion
In light of the preceding two theorems, let us look again at the justification
of `$(p\in\lnot\lnot p)$' (cf. (4.4.10)). First, using (4.3.1), (4.4.5),
(4.4.3), (4.4.40), and (4.3.7), we obtain%
\begin{align*}
(\mathsf{U}  &  \rightarrow\mathsf{J}p,\mathsf{J}\lnot\lnot p\in\mathsf{U}\\
&  \ \ \ \ \ \left.  \wedge%
%TCIMACRO{\tbigwedge }%
%BeginExpansion
{\textstyle\bigwedge}
%EndExpansion
x(x\in\mathsf{J}p\rightarrow p\rightarrow\mathsf{J}\lnot p=0\rightarrow
\lnot\lnot p\rightarrow%
%TCIMACRO{\tbigvee }%
%BeginExpansion
{\textstyle\bigvee}
%EndExpansion
y(y\in\mathsf{J}\lnot\lnot p))\right. \\
&  \rightarrow%
%TCIMACRO{\tbigvee }%
%BeginExpansion
{\textstyle\bigvee}
%EndExpansion
f\in\mathsf{On\ J}p\cap\mathsf{To\ J}\lnot\lnot p%
%TCIMACRO{\tbigwedge }%
%BeginExpansion
{\textstyle\bigwedge}
%EndExpansion
x\in\mathsf{J}p(.fx\in\mathsf{J}\lnot\lnot p)\\
&  \rightarrow%
%TCIMACRO{\tbigvee }%
%BeginExpansion
{\textstyle\bigvee}
%EndExpansion
f(\mathsf{on\ J}p\ \mathsf{to\ J}\lnot\lnot p\ \mathsf{is\ }f)\\
&  \rightarrow%
%TCIMACRO{\tbigvee }%
%BeginExpansion
{\textstyle\bigvee}
%EndExpansion
f\,\mathsf{jst\ }f(p\rightarrow\lnot\lnot p)).
\end{align*}
However, this argument does not produce an explicit description of a
justification of `$(p\in\lnot\lnot p)$'. For this, we apply the more explicit
theorem of choice: by (4.3.3), (4.4.39), (4.4.5), and (4.4.9),%
\begin{align*}
&  \left.  (\mathsf{jst\,}F%
%TCIMACRO{\tbigwedge }%
%BeginExpansion
{\textstyle\bigwedge}
%EndExpansion
x\in\mathsf{J}p%
%TCIMACRO{\tbigvee }%
%BeginExpansion
{\textstyle\bigvee}
%EndExpansion
y(y\in\mathsf{J}\lnot\lnot p)\wedge\phi=\rightthreetimes x\in\mathsf{J}%
p\ \mathsf{The\ }P\ \mathsf{jst\ }P(x\in\mathsf{J}p)\right. \\
&  \ \ \ \ \left.  \rightarrow\mathsf{J}\lnot\lnot p\in\mathsf{U}\right. \\
&  \ \ \ \ \left.  \rightarrow(f=\rightthreetimes x\in\mathsf{J}%
p\ \mathsf{crd}^{\prime}.F(x,\mathsf{.}\phi x)\right. \\
&  \ \ \ \ \ \ \ \ \ \ \ \ \ \ \left.  \rightarrow f\in\mathsf{On\ J}p\wedge%
%TCIMACRO{\tbigwedge }%
%BeginExpansion
{\textstyle\bigwedge}
%EndExpansion
x\in\mathsf{J}p(.fx\in\mathsf{J}\lnot\lnot p)\right. \\
&  \ \ \ \ \ \ \ \ \ \ \ \ \ \ \ \ \ \ \left.  \wedge%
%TCIMACRO{\tbigwedge }%
%BeginExpansion
{\textstyle\bigwedge}
%EndExpansion
x(x\in\mathsf{J}p\rightarrow\lnot\lnot p\rightarrow\mathsf{J}\lnot\lnot
p=1)\right. \\
&  \ \ \ \ \ \ \ \ \ \ \ \ \ \ \left.  \rightarrow f\in\mathsf{On\ J}p\wedge%
%TCIMACRO{\tbigwedge }%
%BeginExpansion
{\textstyle\bigwedge}
%EndExpansion
x\in\mathsf{J}p(.fx\in1)\right. \\
&  \ \ \ \ \ \ \ \ \ \ \ \ \ \ \left.  \rightarrow f\in\mathsf{On\ J}%
p\cap\mathsf{To\ }1\wedge f\in\mathsf{On\ J}p\wedge%
%TCIMACRO{\tbigwedge }%
%BeginExpansion
{\textstyle\bigwedge}
%EndExpansion
x\in\mathsf{J}p(.fx=0)\right. \\
&  \ \ \ \ \ \ \ \ \ \ \ \ \ \ \left.  \rightarrow f=\rightthreetimes
x\in\mathsf{J}p0).\right.
\end{align*}
%

%TCIMACRO{\TeXButton{medskip}{\medskip}}%
%BeginExpansion
\medskip
%EndExpansion%
%TCIMACRO{\TeXButton{noindent}{\noindent}}%
%BeginExpansion
\noindent
%EndExpansion
Comparing this with (4.4.10), we are pleased to have arrived at the same
characterisation of justifications of `$(p\rightarrow\lnot\lnot p)$'$\ $with
or without the exercise of axiom (4.3.10). It may be that excluding (4.3.10)
from our list of axioms would lose us little or nothing other than the
principles of choice (4.4.39) and (4.4.40). However, the axiom does capture a
property of justification sets that merits further exploration.

\subsection{Conclusions}%

%TCIMACRO{\TeXButton{noindent}{\noindent}}%
%BeginExpansion
\noindent
%EndExpansion
We now make the promised return to the question of why axioms (4.3.8) and
(4.3.9) deal only with restricted quantification over a set belonging to the
universe. For example, why don't we replace (4.3.8) by%
\begin{equation}
(\mathsf{jst\,}F%
%TCIMACRO{\tbigwedge }%
%BeginExpansion
{\textstyle\bigwedge}
%EndExpansion
x\in\mathsf{U}\,\underline{\mathsf{u}}x\leftrightarrow\mathsf{on\ dd\smallskip
\,U\ is\ }F\wedge%
%TCIMACRO{\tbigwedge }%
%BeginExpansion
{\textstyle\bigwedge}
%EndExpansion
t\in\mathsf{dd\,U}\ \mathsf{jst\,}.Ft\underline{\mathsf{u}}\mathsf{crd}%
^{\prime}t)? \tag{1}\label{1}%
\end{equation}
If we did so, then, referring to (4.3.1), (4.3.2), (2.20.64), and (2.13.13),
we would obtain%
\begin{align*}
(\mathsf{U}  &  \rightarrow%
%TCIMACRO{\tbigwedge }%
%BeginExpansion
{\textstyle\bigwedge}
%EndExpansion
x\in\mathsf{U}(x=x)\\
&  \rightarrow%
%TCIMACRO{\tbigvee }%
%BeginExpansion
{\textstyle\bigvee}
%EndExpansion
F\in\mathsf{U}\ \mathsf{jst\,}F%
%TCIMACRO{\tbigwedge }%
%BeginExpansion
{\textstyle\bigwedge}
%EndExpansion
x\in\mathsf{U}\,(x=x)\\
&  \rightarrow%
%TCIMACRO{\tbigvee }%
%BeginExpansion
{\textstyle\bigvee}
%EndExpansion
F(F\in\mathsf{U}\wedge\mathsf{on\ U\ is\ }F\wedge%
%TCIMACRO{\tbigwedge }%
%BeginExpansion
{\textstyle\bigwedge}
%EndExpansion
t\in\mathsf{dd\,U}\ \mathsf{jst\,}.Ft(\mathsf{crd}^{\prime}t=\mathsf{crd}%
^{\prime}t)\\
&  \rightarrow%
%TCIMACRO{\tbigvee }%
%BeginExpansion
{\textstyle\bigvee}
%EndExpansion
F(\mathsf{function\ is\ }F\wedge\mathsf{dmn\ }F=\mathsf{U}\wedge
F\in\mathsf{U})\\
&  \rightarrow%
%TCIMACRO{\tbigvee }%
%BeginExpansion
{\textstyle\bigvee}
%EndExpansion
F(\mathsf{function\ is\ }F\wedge\mathsf{U}=\mathsf{dmn\ }F\in\mathsf{U})\\
&  \rightarrow\mathsf{U}\in\mathsf{U}\\
&  \rightarrow0).
\end{align*}
%

%TCIMACRO{\TeXButton{smallskip}{\smallskip}}%
%BeginExpansion
\smallskip
%EndExpansion%
%TCIMACRO{\TeXButton{noindent}{\noindent}}%
%BeginExpansion
\noindent
%EndExpansion
On the other hand, if we adopt%
\begin{align}
&  \left.  (\mathsf{jst\,}T%
%TCIMACRO{\tbigvee }%
%BeginExpansion
{\textstyle\bigvee}
%EndExpansion
x\in A\,\underline{\mathsf{u}}x\right. \nonumber\\
&  \ \ \ \ \left.  \leftrightarrow\mathsf{orderedtriple\ is\ }T\wedge
\mathsf{jst\,crd}^{\prime\prime}T(\mathsf{crd}^{\prime}T\in A)\wedge
\mathsf{jst\,crd}^{\prime\prime\prime}T\underline{\mathsf{u}}\mathsf{crd}%
^{\prime}T)\right.  \tag{2}\label{2}%
\end{align}
instead of (4.3.9), then, using (4.3.1), (2.16.79), (4.4.1), and (2.12.7), we
have
\begin{align*}%
%TCIMACRO{\tbigwedge }%
%BeginExpansion
{\textstyle\bigwedge}
%EndExpansion
x(x\in\mathsf{U}  &  \rightarrow x\in\mathsf{U}\wedge x=x\\
&  \rightarrow%
%TCIMACRO{\tbigvee }%
%BeginExpansion
{\textstyle\bigvee}
%EndExpansion
\alpha,\beta(\alpha\in\mathsf{J}(x\in\mathsf{U})\wedge\beta\in\mathsf{J}%
(x=x))\\
&  \rightarrow%
%TCIMACRO{\tbigvee }%
%BeginExpansion
{\textstyle\bigvee}
%EndExpansion
\alpha,\beta(T=x,\alpha,\beta\rightarrow\mathsf{orderedtriple\ is\ }%
T\wedge\mathsf{crd}^{\prime\prime}T\in\mathsf{J}(\mathsf{crd}^{\prime}%
T\in\mathsf{U})\\
&  \left.  \wedge\ \mathsf{crd}^{\prime\prime\prime}T\in\mathsf{J}%
(\mathsf{crd}^{\prime}T=\mathsf{crd}^{\prime}T))\right. \\
&  \rightarrow%
%TCIMACRO{\tbigvee }%
%BeginExpansion
{\textstyle\bigvee}
%EndExpansion
\alpha,\beta%
%TCIMACRO{\tbigvee }%
%BeginExpansion
{\textstyle\bigvee}
%EndExpansion
T(\mathsf{orderedtriple\ is\ }T\wedge\mathsf{crd}^{\prime}T=x\in
\mathsf{U}\wedge\mathsf{crd}^{\prime\prime}T\in\mathsf{J}(\mathsf{crd}%
^{\prime}T\in\mathsf{U})\\
&  \left.  \wedge\ \mathsf{crd}^{\prime\prime\prime}T\in\mathsf{J}%
(\mathsf{crd}^{\prime}T=\mathsf{crd}^{\prime}T))\right. \\
&  \rightarrow%
%TCIMACRO{\tbigvee }%
%BeginExpansion
{\textstyle\bigvee}
%EndExpansion
T(\mathsf{orderedtriple\ is\ }T\wedge T\in\mathsf{U}\wedge\mathsf{crd}%
^{\prime}T=x\wedge\mathsf{crd}^{\prime\prime}T\in\mathsf{J}(\mathsf{crd}%
^{\prime}T\in\mathsf{U})\\
&  \left.  \wedge\ \mathsf{crd}^{\prime\prime\prime}T\in\mathsf{J}%
(\mathsf{crd}^{\prime}T=\mathsf{crd}^{\prime}T))\right. \\
&  \rightarrow%
%TCIMACRO{\tbigvee }%
%BeginExpansion
{\textstyle\bigvee}
%EndExpansion
T(\mathsf{jst\ }T%
%TCIMACRO{\tbigvee }%
%BeginExpansion
{\textstyle\bigvee}
%EndExpansion
z\in\mathsf{U}(z=z)\wedge\mathsf{crd}^{\prime}T=x)\\
&  \rightarrow x\in\mathsf{U}\wedge%
%TCIMACRO{\tbigvee }%
%BeginExpansion
{\textstyle\bigvee}
%EndExpansion
T\in\mathsf{J}%
%TCIMACRO{\tbigvee }%
%BeginExpansion
{\textstyle\bigvee}
%EndExpansion
z\in\mathsf{U}(z=z)(\mathsf{crd}^{\prime}T=x)\\
&  \rightarrow x\in\mathsf{E}x%
%TCIMACRO{\tbigvee }%
%BeginExpansion
{\textstyle\bigvee}
%EndExpansion
T\in\mathsf{J}%
%TCIMACRO{\tbigvee }%
%BeginExpansion
{\textstyle\bigvee}
%EndExpansion
z\in\mathsf{U}(z=z)(\mathsf{crd}^{\prime}T=x)).
\end{align*}
Hence%
\[
(\mathsf{U\subset E}x%
%TCIMACRO{\tbigvee }%
%BeginExpansion
{\textstyle\bigvee}
%EndExpansion
T\in\mathsf{J}%
%TCIMACRO{\tbigvee }%
%BeginExpansion
{\textstyle\bigvee}
%EndExpansion
z\in\mathsf{U}(z=z)(\mathsf{crd}^{\prime}T=x)\subset\mathsf{U})
\]
and therefore%
\[
(\mathsf{U}=\mathsf{E}x%
%TCIMACRO{\tbigvee }%
%BeginExpansion
{\textstyle\bigvee}
%EndExpansion
T\in\mathsf{J}%
%TCIMACRO{\tbigvee }%
%BeginExpansion
{\textstyle\bigvee}
%EndExpansion
z\in\mathsf{U}(z=z)(\mathsf{crd}^{\prime}T=x)).
\]
Now,%
\[%
%TCIMACRO{\tbigwedge }%
%BeginExpansion
{\textstyle\bigwedge}
%EndExpansion
T\in\mathsf{J}%
%TCIMACRO{\tbigvee }%
%BeginExpansion
{\textstyle\bigvee}
%EndExpansion
z\in\mathsf{U}(z=z)(\mathsf{crd}^{\prime}T\in\mathsf{U}),
\]
so, by (2.20.186), (4.3.3), (2.20.63), and (2.13.13),%
\begin{align*}
(f  &  =\rightthreetimes T\in\mathsf{J}%
%TCIMACRO{\tbigvee }%
%BeginExpansion
{\textstyle\bigvee}
%EndExpansion
z\in\mathsf{U}(z=z)\mathsf{crd}^{\prime}T\\
&  \rightarrow\mathsf{function\ is\ }f\\
&  \left.  \wedge\ \mathsf{dmn\ }f=\mathsf{E}T\in\mathsf{J}%
%TCIMACRO{\tbigvee }%
%BeginExpansion
{\textstyle\bigvee}
%EndExpansion
z\in\mathsf{U}(z=z)(\mathsf{crd}^{\prime}T\in\mathsf{U})\right. \\
&  \left.  \wedge\
%TCIMACRO{\tbigwedge }%
%BeginExpansion
{\textstyle\bigwedge}
%EndExpansion
T\in\mathsf{dmn\ }f(.fx=\mathsf{crd}^{\prime}T)\right. \\
&  \rightarrow\mathsf{function\ is\ }f\wedge\mathsf{dmn\ }f=\mathsf{J}%
%TCIMACRO{\tbigvee }%
%BeginExpansion
{\textstyle\bigvee}
%EndExpansion
z\in\mathsf{U}(z=z)\in\mathsf{U}\\
&  \left.  \wedge\ \mathsf{rng\ }f=\mathsf{E}x%
%TCIMACRO{\tbigvee }%
%BeginExpansion
{\textstyle\bigvee}
%EndExpansion
T\in\mathsf{J}%
%TCIMACRO{\tbigvee }%
%BeginExpansion
{\textstyle\bigvee}
%EndExpansion
z\in\mathsf{U}(z=z)(\mathsf{crd}^{\prime}T=x)\right. \\
&  \rightarrow\mathsf{function\ is\ }f\wedge\mathsf{dmn\ }f\in\mathsf{U}%
\wedge\mathsf{rng\ }f=\mathsf{U}\\
&  \rightarrow\mathsf{rng\ }f\in\mathsf{U\wedge rng\ }f=\mathsf{U}\\
&  \rightarrow\mathsf{U}\in\mathsf{U}\\
&  \rightarrow0).
\end{align*}
We conclude that each of (\ref{1}) and (\ref{2}) leads to a contradiction, so
neither can be used as an axiom.%

%TCIMACRO{\TeXButton{smallskip}{\smallskip}}%
%BeginExpansion
\smallskip
%EndExpansion%
%TCIMACRO{\TeXButton{noindent}{\noindent}}%
%BeginExpansion
\noindent
%EndExpansion
In view of this apparent drawback to our axiomatic constructive justification
theory, it is worth pointing out that in the practice of constructive
mathematics \`{a} la Bishop \cite{Bishop,BB}, the domain of a quantification
is generally some well-constructed set (that is, a member of the universe).
Even in CMST, a formal foundation for Bishop-style constructive mathematics,
unrestricted quantifications are often of a trivial nature, like `$%
%TCIMACRO{\tbigwedge }%
%BeginExpansion
{\textstyle\bigwedge}
%EndExpansion
x(x=x)$', and can be replaced by free-variable statements (like `$(x=x)$'). An
obvious exception to this is the definition `$(0\equiv%
%TCIMACRO{\tbigwedge }%
%BeginExpansion
{\textstyle\bigwedge}
%EndExpansion
xx)$'. One way out of this would be to take `$0$' as a primitive constant and
make `$(%
%TCIMACRO{\tbigwedge }%
%BeginExpansion
{\textstyle\bigwedge}
%EndExpansion
xx\equiv0)$' an axiom. But we shall content ourselves with the situation
described in \cite{BAlps,dsb2022} and this paper; so our justification theory
does not incorporate unrestricted universal quantification.%

%TCIMACRO{\TeXButton{newpage}{\newpage}}%
%BeginExpansion
\newpage
%EndExpansion

\section{Appendix: \ Axioms for Logic}%

%TCIMACRO{\TeXButton{medskip}{\medskip}}%
%BeginExpansion
\medskip
%EndExpansion
%

%TCIMACRO{\TeXButton{noindent}{\noindent}}%
%BeginExpansion
\noindent
%EndExpansion
\textbf{Axioms for propositional logic\ }

\begin{enumerate}
\item $\left(  q\rightarrow\left(  p\rightarrow q\right)  \right)  $

\item $\left(  \left(  p\rightarrow\left(  q\rightarrow r\right)  \right)
\rightarrow\left(  \left(  p\rightarrow q\right)  \rightarrow\left(
p\rightarrow r\right)  \right)  \right)  $

\item $\left(  \left(  p\wedge q\right)  \rightarrow p\right)  $

\item $\left(  \left(  p\wedge q\right)  \rightarrow q\right)  $

\item $\left(  p\rightarrow\left(  q\rightarrow(p\wedge q)\right)  \right)  $

\item $\left(  p\rightarrow\left(  p\vee q\right)  \right)  $

\item $\left(  q\rightarrow\left(  p\vee q\right)  \right)  $

\item $\left(  \left(  \left(  p\rightarrow r\right)  \wedge\left(
q\rightarrow r\right)  \right)  \rightarrow\left(  \left(  p\vee q\right)
\rightarrow r\right)  \right)  $
\end{enumerate}

%

%TCIMACRO{\TeXButton{bigskip}{\bigskip}}%
%BeginExpansion
\bigskip
%EndExpansion
%

%TCIMACRO{\TeXButton{noindent}{\noindent}}%
%BeginExpansion
\noindent
%EndExpansion
\textbf{Axioms for predicate logic}

The predicate-logic axioms consist of the propositional axioms
(1.4.1)--(1.4.8) and the following:%
\index{Axiom(s)!for predicate logic}%
%TCIMACRO{\TeXButton{smallskip}{\smallskip}}%
%BeginExpansion
\smallskip
%EndExpansion

\begin{enumerate}
\item $\left(
%TCIMACRO{\tbigwedge }%
%BeginExpansion
{\textstyle\bigwedge}
%EndExpansion
x\,\underline{\mathsf{u}}x\rightarrow\underline{\mathsf{u}}x\right)  $

\item $\left(  \underline{\mathsf{u}}x\rightarrow%
%TCIMACRO{\tbigvee }%
%BeginExpansion
{\textstyle\bigvee}
%EndExpansion
x\,\underline{\mathsf{u}}x\right)  $

\item $\left(  y\rightarrow%
%TCIMACRO{\tbigwedge }%
%BeginExpansion
{\textstyle\bigwedge}
%EndExpansion
xy\right)  $

\item $\left(
%TCIMACRO{\tbigvee }%
%BeginExpansion
{\textstyle\bigvee}
%EndExpansion
xy\rightarrow y\right)  $

\item $\left(
%TCIMACRO{\tbigwedge }%
%BeginExpansion
{\textstyle\bigwedge}
%EndExpansion
x\left(  \underline{\mathsf{u}}x\rightarrow\underline{\mathsf{v}}x\right)
\rightarrow\left(
%TCIMACRO{\tbigwedge }%
%BeginExpansion
{\textstyle\bigwedge}
%EndExpansion
x\,\underline{\mathsf{u}}x\rightarrow%
%TCIMACRO{\tbigwedge }%
%BeginExpansion
{\textstyle\bigwedge}
%EndExpansion
x\,\underline{\mathsf{v}}x\right)  \right)  $

\item $(%
%TCIMACRO{\tbigwedge }%
%BeginExpansion
{\textstyle\bigwedge}
%EndExpansion
x\left(  \underline{\mathsf{u}}x\rightarrow\underline{\mathsf{v}}x\right)
\rightarrow\left(
%TCIMACRO{\tbigvee }%
%BeginExpansion
{\textstyle\bigvee}
%EndExpansion
x\,\underline{\mathsf{u}}x\rightarrow%
%TCIMACRO{\tbigvee }%
%BeginExpansion
{\textstyle\bigvee}
%EndExpansion
x\,\underline{\mathsf{v}}x\right)  )$
\end{enumerate}

%

%TCIMACRO{\TeXButton{newpage}{\newpage}}%
%BeginExpansion
\newpage
%EndExpansion

%

%TCIMACRO{\TeXButton{bigskip}{\bigskip}}%
%BeginExpansion
\bigskip
%EndExpansion
%

%TCIMACRO{\TeXButton{noindent}{\noindent}}%
%BeginExpansion
\noindent
%EndExpansion
\textbf{Author's address: \ }Department of Mathematics \& Statistics,
University of Canterbury, Christchurch 8140, New Zealand%

%TCIMACRO{\TeXButton{noindent}{\noindent}}%
%BeginExpansion
\noindent
%EndExpansion
\textbf{Author's email: \ \ \ \ \ }\texttt{dugbridges@gmail.com}

\vfill

\begin{flushright}
{\small \texttt{dsbridges 180626}}
\end{flushright}


\begin{thebibliography}{99}                                                                                               %


\bibitem {Alps}\textsf{R.A. Alps (ed.), \textit{A.P. Morse's Set Theory and
Analysis}, Birkh\"{a}user, Cham, Switzerland, 2022. }

\bibitem {BAlps}R.A. Alps and D.S. Bridges, \textit{Constructive Morse Set
Theory---A Foundation for Constructive Mathematics}, monograph, in preparation.

\bibitem {Beeson}\textsf{M. Beeson, \emph{Foundations of Constructive
Mathematics}, Springer Verlag, Heidelberg, 1985. }

\bibitem {Bishop}E. Bishop, \emph{Foundations of Constructive Analysis,}
McGraw-Hill, New York, 1967.

\bibitem {Bishnum}E. Bishop, `Mathematics as a Numerical Language', in
\textit{Intuitionism and Proof Theory} (A. Kino, J. Myhill, R.E. Vesley, eds),
North-Holland Pub. Co., Amsterdam, 1970.

\bibitem {BB}E. Bishop and D.S. Bridges, \emph{Constructive Analysis},
Grundlehren der Math. Wiss. \textbf{279}, Springer Verlag, Heidelberg, 1985.

\bibitem {dsb2022}D.S.\ Bridges, `Morse set theory as a foundation for
constructive mathematics', Theoretical Comp. Sci. 928, 115--135, 2022. https://doi.org/10.1016/j.tcs.2022.06.019

\bibitem {Diac}\textsf{R. Diaconescu, `Axiom of choice and complementation',
Proc. Amer. Math. Soc. \textbf{51}, 176--178, 1975.}

\bibitem {Goodmy}\textsf{N. Goodman and J. Myhill, `Choice implies excluded
middle', Zeit. Math. Logik Grundlagen Math. \textbf{23}, 461, 1978.}

\bibitem {KV}S..C. Kleene and R.E. Vesley, \emph{The Foundations of
Intuitionistic Mathematics, Especially in Relation to Recursive Functions,
}North-Holland, Amsterdam, 1965.

\bibitem {Morse}A.P. Morse, \emph{A\ Theory of Sets}, Academic Press, New
York, 1965 (second edition 1986).

\bibitem {Morse2}\textsf{A.P. Morse, \emph{A Theory of Sets} (second edition),
Academic Press, New York, 1986.}

\bibitem {TvD}A.S. Troelstra and D. van Dalen, \emph{Constructivism in
Mathematics} (Vols 1 and 2), North Holland Publ. Co., Amsterdam, 1988.
\end{thebibliography}
\end{document}